\documentclass{amsart}
\usepackage{amssymb,color}
\usepackage{mathrsfs,epsfig}
\date{October 10, 2001}
\title[Spectral Theorem for Convex Monotone Homogeneous maps]{Spectral Theorem for Convex Monotone Homogeneous Maps, and Ergodic Control}
\thanks{This work was partially supported by the European Community Framework IV program through the research network ALAPEDES (``The Algebraic Approach to Performance Evaluation of Discrete Event Systems'').
The results of this paper
were announced in~\cite{gaubert01b}.}
\author{Marianne Akian}
\address{INRIA, Domaine de Voluceau, 78153 Le Chesnay C\'edex, France}
\email{Marianne.Akian@inria.fr}
\author{St\'ephane Gaubert}
\address{ENSTA, 32 Bd. Victor, 78539 Paris C\'edex 15,  France,
and 
INRIA, Domaine de Voluceau, 78153 Le Chesnay C\'edex, France.}
\email{Stephane.Gaubert@inria.fr}

\newcommand{\sat}{\mrm{Sat}}
\newcommand{\ov}[1]{\overline{#1}}

\newcommand{\EE}{\mathbb{E}}
\newcommand{\new}[1]{{\em #1}}
\newcommand{\mrm}[1]{\text{\rm #1}}
\newcommand{\sfr}{\mathsf{r}}
\newcommand{\sfi}{\mathsf{i}}
\newcommand{\vex}{\operatorname{co}}
\newcommand{\vexb}{\operatorname{\overline{co}}}
\newcommand{\lcm}{\text{\rm lcm}}
\newcommand{\sfm}{\mathsf{m}}
\newcommand{\sfc}{\mathsf{c}}
\newcommand{\dom}{\operatorname{dom}}
\newcommand{\cl}{\operatorname{cl}}
\newcommand{\access}{\stackrel{*}{\rightarrow}}
\newcommand{\bydef}{\stackrel{\text{\rm def}}{=}}
\newcommand{\R}{\mathbb{R}}
\newcommand{\Rmin}{\R\cup\{+\infty\}}
\newcommand{\Rmax}{\R\cup\{-\infty\}}
\newcommand{\N}{\mathbb{N}}
\newcommand{\Rp}{\mathbb{R}_{^{^+}\!\!}^*}
\newcommand{\sS}{\mathscr{S}}
\newcommand{\proba}[1]{\sS_{\!#1}}
\newcommand{\probam}[1]{\sS_{\!#1}^{-}}
\newcommand{\sE}{\mathscr{E}}
\newcommand{\sN}{\mathscr{N}}
\newcommand{\sF}{\mathscr{F}}
\newcommand{\sC}{\mathscr{C}}
\newcommand{\sG}{\mathscr{G}}
\newcommand{\sL}{\mathscr{L}}
\newcommand{\sM}{\mathscr{M}}
\newcommand{\sP}{\mathscr{P}}
\newcommand{\sQ}{\mathscr{Q}}
\newcommand{\sR}{\mathscr{R}}
\newcommand{\nf}{\mathsf{N}^{\rm f}}
\newcommand{\nc}{\mathsf{N}^{\rm c}}
\newcommand{\gf}{\sG^{\rm f}}
\newcommand{\gc}{\sG^{\rm c}}
\newcommand{\cf}{\sC^{\rm f}}
\newcommand{\cc}{\sC^{\rm c}}
\newcommand{\ec}{\sE^{\rm c}}
\newcommand{\set}[2]{\{#1\mid\,#2\}}
\newcommand{\comp}{\circ}
\newtheorem{theorem}{Theorem}[section]
\newtheorem{prop}[theorem]{Proposition}
\newtheorem{lemma}[theorem]{Lemma}
\newtheorem{cor}[theorem]{Corollary}
\theoremstyle{remark}
\newtheorem{example}[theorem]{Example}
\theoremstyle{definition}
\newtheorem{defi}[theorem]{Definition}
\begin{document}
\sloppy
\begin{abstract}
We consider
convex maps $f:\R^n\to\R^n$ 
that are monotone
(i.e., that preserve
the product ordering of $\R^n$), and 
nonexpansive for the sup-norm.
This includes convex monotone maps that are
additively homogeneous (i.e., that commute
with the addition of constants).
We show that the fixed point set of $f$,
when it is non-empty, is isomorphic to a convex inf-subsemilattice
of $\R^n$, whose dimension is at most equal to the number
of strongly connected components of a critical
graph defined from the tangent affine maps of $f$.
This yields in particular an uniqueness result
for the bias vector of ergodic control problems.
This generalizes results obtained previously by Lanery,
Romanovsky, and Schweitzer and Federgruen, for ergodic control problems with
finite state and action spaces, which correspond to the special case
of piecewise affine maps $f$. 
We also show that the length of periodic orbits
of $f$ is bounded by the cyclicity of its
critical graph, which implies that the possible orbit
lengths of $f$ are exactly the orders of elements
of the symmetric group on $n$ letters. 
\end{abstract}
\maketitle
\tableofcontents
\section{Introduction}
\subsection{Motivations and Statement of the Main Result}

We say that a map $f:\R^n\to\R^n$ is \new{monotone} if for
all $x,y\in \R^n$, $x\leq y
\implies f(x)\leq f(y)$, where $\leq$ denotes the product
ordering of $\R^n$ ($x\leq y$ if $x_i\leq y_i$, for all
$1\leq i\leq n$). We say that $f$ is additively \new{homogeneous}
if for all $\lambda\in \R$,
$x\in \R^n$, $f(\lambda+x)=\lambda +f(x)$, where 
$\lambda+x=(\lambda+x_1,\ldots,\lambda+x_n)$. 
It is easy to see that
a monotone homogeneous map is nonexpansive for the sup-norm:
for all $x,y\in \R^n$, $|f(x)-f(y)|\leq |x-y|$,
where $|x|=\max_{1\leq i\leq n}|x_i|$
(see~\cite{crandall}).

Monotone homogeneous maps arise classically 
in optimal control and game theory (see for instance
 \cite{bewkohl,kohlberg80,whittle,kol92,filar97,sorin}), 
in the modeling of discrete events systems (see 
\cite{BCOQ,guna94,CGQ95a,CGQ95b,vin97,gg98a,mxpnxp0}), 
and in nonlinear potential theory~\cite{dellach96}, as nonlinear extension
of Markov transitions.
They also arise in nonlinear Perron-Frobenius theory, when one considers
multiplicatively homogeneous 
maps $F$ acting on a cone and preserving the order of the cone: in
the simplest case, when the cone is $(\Rp)^n$ 
(where $\Rp=\set{x\in \R}{x>0}$), the transformation
$F\mapsto \log\comp F\comp \exp$
(where $\log:(\Rp)^n\to \R^n$ is the
map which does $\log$ entrywise, and $\exp=\log^{-1}$)
sends the set of monotone multiplicatively
homogeneous maps to the set of monotone
additively homogeneous maps. See
for instance~\cite{bougerol,wojtkowski94,morishima,ms69,nuss88,nuss89,sabot,wojtkowski85}
for various studies and applications. 

A basic problem, for a monotone homogeneous map $f$,
is the existence, and uniqueness (up to an additive
constant), of the \new{additive eigenvectors}
of $f$, which are the $v\in \R^n$ such that $f(v)=\lambda+v$,
for some \new{additive eigenvalue}
$\lambda\in \R$. In the sequel, we will omit
the term ``additive'', when the additive
nature of the objects will be clear 
from the context. 
When $f$ has an eigenvector $v$ with eigenvalue $\lambda$,
by homogeneity of $f$, $f^k(v)=k\lambda +v$ holds for all $k\geq 0$, and
by nonexpansiveness of $f$, $|f^k(x)-k\lambda-v|=|f^k(x)-f^k(v)|
\leq |x-v|$, hence,
\begin{equation}
f^k(x)=k\lambda +O(1) \qquad\text{\rm when $k\to \infty$,}
\label{lambda}
\end{equation}
for all $x\in \R^n$
(all the orbits of $f$
have a linear growth rate of $\lambda$). 
This implies in particular that the eigenvalue $\lambda$
is unique. Hence, we can speak
without ambiguity of the \new{eigenspace} of $f$,
which is the set $\sE(f)=\set{x\in \R^n}{f(x)=\lambda +x}$.
 In many applications, the eigenvalue and eigenvector are 
fundamental objects: for instance, in stochastic control,
the eigenvalue gives the optimal reward per time unit, 
and eigenvectors give stationary rewards (we 
explain this in detail in Section~\ref{sec-stoch}).
In discrete event systems applications, the eigenvalue gives the throughput,
and eigenvectors give stationary schedules. 

Several Perron-Frobenius like theorems
guarantee the existence of eigenvectors
of monotone (additively) homogeneous
maps $\R^n\to\R^n$. Such results go back at least to
Kre{\u\i}n and Rutman~\cite[\S7]{krein},
in the context of monotone multiplicatively
homogeneous maps leaving a cone in a Banach space
invariant, and to Morishima,
whose book~\cite{morishima} contains
a complete study of finite dimensional
non-linear Perron-Frobenius theory. 
A modern overview
appears in the memoirs of Nussbaum~\cite{nuss88,nuss89},
which contain in particular general existence results
for eigenvectors,
following~\cite{nussbaum86}. 
Different existence conditions appeared in~\cite{oshime84}.
General results on the geometry of the eigenspaces are available,
for instance, the result of Bruck~\cite{bruck} shows in particular that
$\sE(f)$ is the image of a nonexpansive projector and a fortiori is 
connected, see also~\cite[Theorems 4.5, 4.6 and 4.7]{nuss88}.

In this paper, we describe the eigenspaces of
\new{convex} monotone homogeneous maps
$f:\R^n\to \R^n$.
(We say that a $\R^n$-valued map is convex
when its coordinates are convex. We refer the reader
to~\cite{ROCK} for all convexity notions used in the paper: subdifferentials,
domain, Fenchel transform, etc.)

To state our main result, we need a few
definitions (see Section~\ref{cs} for details).
We first generalize the notion of subdifferential
to maps $\R^n\to \R^n$ by setting,
for $x\in \R^n$,
$\partial f(x)=\set{P\in \R^{n\times n}}{f(y)-f(x) \geq P(y-x),
\;\forall y\in \R^n}$.
It is easy to see (Corollary~\ref{prop-1} and Equation~\eqref{rect}
below)
that by monotonicity and homogeneity of $f$,
the elements of $\partial f(x)$ are stochastic matrices.
If $v$ is an eigenvector of $f$,
we call \new{critical graph} of $f$,
the (directed) graph $\gc(f)$ which is the union 
of final graphs of stochastic matrices $P\in \partial f(v)$
(we call \new{final graph} of a stochastic
matrix the restriction of its graph to the
set of final classes, see \S\ref{cs2} and~\S\ref{ss-crit}).
The graph $\gc(f)$ is independent
of the choice of the eigenvector $v$
(Proposition~\ref{indep} below).
We call \new{critical nodes} of $f$, the nodes of $\gc(f)$,
and denote by $\nc(f)$ the set of critical nodes.
We call \new{critical classes}
of $f$ the sets of nodes
$C_1,\ldots,C_s$ of
the strongly connected components of $\gc(f)$,
$\sG_1,\ldots,\sG_s$.
We call \new{cyclicity} of $\sG_i$, and denote
by $\sfc(\sG_i)$, the 
$\gcd$ of the lengths of the circuits
of $\sG_i$, and we define the \new{cyclicity} of $f$
by $\sfc(f)=\lcm(\sfc(\sG_1),\ldots,\sfc(\sG_s))$.
We say that a monotone homogeneous map $g:U\subset \R^n\to V\subset\R^p$ is
a monotone homogeneous \new{isomorphism} if it has a monotone
homogeneous inverse.
The following theorem gathers results
from Theorem~\ref{thmain}, Corollary~\ref{cor-summarize},
Corollary~\ref{cor-asymp}, 
and Theorem~\ref{th2} below.
\begin{theorem}[Convex Spectral Theorem]\label{t-cst}
Let $f:\R^n\to\R^n$ denote a convex monotone
homogeneous map that has an eigenvector.
Denote by $C=\nc(f)$ the set of critical nodes of $f$,
$c=\sfc(f)$ the cyclicity of $f$,
and $\lambda$ the unique eigenvalue of $f$.
Then, 
\begin{enumerate}
\item the restriction $\sfr:\R^n\to \R^C,
x\mapsto (x_i)_{i\in C}$,
is a monotone homogeneous isomorphism
from $\sE(f)$ to its image $\ec(f)$;
\item $\ec(f)$ is an inf-subsemilattice of $(\R^C,\leq)$;
\item $\ec(f)$ is a convex set whose dimension is at most
equal to the number of critical classes of
$f$, and this bound is attained 
when $f$ is piecewise affine;
\item for all $x\in \R^n$,
$f^{kc}(x)-kc\lambda$ has a limit when $k\to\infty$.
\end{enumerate}
\end{theorem}
In particular, when $f$ has only one critical
class, the eigenvector of $f$ is unique (up to an additive constant).
It also follows from the last
assertion of the theorem that the set of limit points of 
$f^{k}(x)-k\lambda$ when $k\to \infty$, and
$x\in \R^n$, is precisely $\sE(f^c)$.
Theorem~\ref{t-cst} also allows
us to bound the dimension of this set.
Indeed, we shall see in Theorem~\ref{ncf} and Proposition~\ref{orbitclass}
below that the set of critical nodes is the same
for $f$ and $f^c$, and that $f^c$ has
$\sfc(\sG_1)+ \cdots +\sfc(\sG_s)$
critical classes. Hence, applying Theorem~\ref{t-cst}
to $f^c$, we get that the restriction $\sfr$ 
is a monotone homogeneous isomorphism 
from $\sE(f^c)$ to a convex set, $\ec(f^c)$,
of dimension at most $\sfc(\sG_1)+ \cdots +\sfc(\sG_s)$,
the bound being attained when $f$ is piecewise
affine. 

The paper is devoted to the proof (Sections~2--6) and to the
stochastic control interpretation (Section~7) of Theorem~\ref{t-cst}.
In Section~\ref{cs}, we detail the definitions and properties of 
subdifferentials  and critical graph  of
convex monotone homogeneous maps.
An important element of the proofs is 
the maximum principle for Markov chains
(Lemma~\ref{disc}). 
In Section~\ref{eigen}, we establish the first part of Theorem~\ref{t-cst}
concerning the structure of the eigenspace:
points 1, 2 and the first assertion in point 3.
The main argument is again the maximum principle.
Section~\ref{secttang} is devoted to further tools and properties 
used in the remaining sections, which are of independent interest:
Theorem~\ref{ncf} shows that $\gc(f^k)=\gc(f)^k$
(this will be used in Section~\ref{cyclicity} for the proof of the
cyclicity part  of Theorem~\ref{t-cst}; we also introduce directional 
derivatives (which will be used in Section~\ref{piecewise} to
connect  $\sE(f)$ to  $\sE(f'_v)$ for any eigenvector
$v$), additive recession functions 
(formula~\eqref{defghat2}),
invariant critical classes and the associated
decomposition of $f$ (Lemma~\ref{lem-restriction}),
and also a characterization of the set
of critical nodes in terms of supports
of nonlinear ``excessive'' measures
(Proposition~\ref{mesf}).

Section~\ref{cyclicity} is devoted to the proof of point
4 of Theorem~\ref{t-cst}.
This result relies on a more general theorem
of Nussbaum~\cite{nuss90} and Sine~\cite{sine},
which states that if $f:\R^n\to\R^n$ is nonexpansive for the sup-norm
and has a fixed point, then, for all $x\in \R^n$,
$f^{kc}(x)$ converges when $k\to \infty$,
for some minimal constant $c$ which can be bounded
by a function of $n$. When $f$ is convex 
monotone and homogeneous,
the last assertion of the convex spectral
theorem shows that the possible values of $c$
are exactly the orders of elements of the symmetric group
on $n$ letters. Equivalently,
convex monotone homogeneous maps 
have the same orbit lengths as permutation matrices,
 a result which was known to be true in the special cases of linear maps
associated to nonnegative matrices (see~\cite[chapter 9]{nuss99}),
of linear maps
over the max-plus semiring (see~\cite{cohen83,nussbaum91}),
and also of piecewise affine convex
monotone homogeneous maps (which include max-plus linear maps),
see the discussion in \S\ref{sec-rel}
below. More generally, computing the orbit lengths of nonexpansive maps
for polyhedral norms raises interesting combinatorial and analytical
problems (see in 
particular~\cite{Akcoglu,weller,scheu88,scheu91,nuss90,nuss91b,nuss98,nuss99,lemmens}).

The equality in point 3 of Theorem~\ref{t-cst} is proved in  
Section~\ref{piecewise}. As will be  discussed in 
\S\ref{sec-rel}, this part of the theorem has already
been proved by Romanovsky~\cite{roma73}
and by Schweitzer and Federgruen~\cite{sch78}.
We provide here an independent proof, which emphasizes the 
qualitative properties of $\sE(f)$, using the tools of 
Section~\ref{secttang}. We also give  a  polynomial
time algorithm to compute $\gc(f)$.

\subsection{Related Optimal Control Results}\label{sec-rel}
Convex monotone homogeneous maps $f:\R^n\to\R^n$ are exactly 
dynamic programming operators associated to stochastic control
problems with state space $\{1,\ldots,n\}$. Computing the stationary solutions
and the asymptotic behavior of solutions of 
dynamic programming equations is an old
problem of stochastic control which is essentially equivalent to that
of computing the eigenspace $\sE(f)$ and the asymptotics of
$f^k$ when $k\to\infty$.
This has been much studied in the stochastic control literature,
particularly in the case of finite action spaces, which corresponds
to piecewise affine maps.
In this special case, results equivalent to the third assertion of 
Theorem~\ref{t-cst} were obtained by Romanovsky~\cite{roma73}
using linear programming techniques, and also by Schweitzer and 
Federgruen~\cite{sch78} who gave an explicit representation of $\sE(f)$
in terms of resolvents associated to optimal strategies (see 
\cite[Theorem 4.1]{sch78}).
Again in this special case, a result equivalent to the fourth assertion
of Theorem~\ref{t-cst}  was stated by Lanery~\cite{lanery}.
The arguments 
of~\cite{lanery} only proved the special case where $\nc(f)=\{1,\ldots,n\}$,
see the discussion in the introduction and in Note~1 of~\cite{sch77}. 
A proof valid for a general $\nc(f)\subset \{1,\ldots,n\}$
was given by Schweitzer and Federgruen~\cite{sch77},
who also proved the optimality of $\sfc(f)$.

The special case of deterministic control problems leads to maps $f$
that are max-plus linear. These maps have been studied independently
by the max-plus community.
In this context, the dimension of the eigenspace was characterized
by Gondran and Minoux~\cite{gondran77},
and the remaining part of the max-plus spectral theorem, dealing
with cyclicity, was obtained by
Cohen, Dubois, Quadrat, and Viot~\cite{cohen83}
(see also~\cite{BCOQ}).
(Note however that more precise results --explicit form of the
eigenspace, finite time convergence of the iterates-- are available
in the max-plus case.)
The max-plus spectral theorem has a long story, which goes back to
Cuninghame-Green (see~\cite{cuning79}
and the references therein), Romanovsky~\cite{romanovski},
and Vorobyev~\cite{vorobyev67},
to quote the most ancient contributions.
See the collection of articles~\cite{maslov92},
Kolokoltsov and Maslov~\cite{maslovkololtsov95},
and the references therein,
for generalizations to infinite dimension.
See also~\cite{GM,bapat98} for surveys.

The present work was inspired by the max-plus spectral theorem
and uses nonexpansive maps techniques (we were
unaware of the results of~\cite{lanery,roma73,sch77,sch78}).
We next emphasize differences with earlier results.
We consider general convex
monotone homogeneous maps, which
correspond to stochastic control problems
with finite state space and {\em arbitrary}
action spaces, whereas the results in~\cite{lanery,roma73,sch77,sch78}
require the action space to be finite.
Our proof technique, which relies on the maximum principle,
 can be naturally transposed to other (infinite dimensional) contexts.
Another tool in our proof is
the critical graph $\gc(f)$, which generalizes the critical graph that
appears in max-plus algebra (see Proposition~\ref{max-crit}).
The critical graph already
appeared in~\cite[p. 491]{roma73}, with a different definition
in terms of optimal policies. The new definition that we give here in terms
of subdifferentials leads in particular to a polynomial
time algorithm to compute $\gc(f)$ (see~\S\ref{sect-algo}).
(The equivalence of both definitions is shown
in Proposition~\ref{t-o}.)  It should
also be noted that when passing from
the case of finite action spaces
to arbitrary action spaces, new
phenomena occur. For instance, 
Example~\ref{moreau} shows that
we cannot hope, in general, to
characterize the dimension of eigenspaces
in terms of graphs like $\gc(f)$.

Let us mention in passing that the critical graph
has an intuitive interpretation in terms of ``recurrence''.
For a Markov chain, a node is recurrent if the probability of return
to this node is equal to one.
For a max-plus matrix with eigenvalue $0$,
a node is ``recurrent'', i.e. belong to a critical class,
if we can return to this node with zero reward.
When $f$ is a convex monotone homogeneous map with eigenvalue $0$,
a node $i$ is ``recurrent'', i.e. belong to a critical class,
if we can find a strategy for which, starting from $i$,
we eventually return to $i$ with probability $1$ and zero
mean reward. This provides a new illustration
of the analogy between probability and optimization
developed in~\cite[Chap.~VIII]{maslov73}, \cite{salut90,jpq90,aqv95,delmoral95}, \cite[\S4.2]{jpqmaxplus97}, \cite{litvinov00}, and~\cite{tolya01}. 

Ergodic control problems of diffusion processes lead to spectral problems for
infinite dimensional monotone homogeneous semigroups which can be expressed
in terms of ergodic Hamilton-Jacobi-Bellman (HJB) 
partial differential equations.
In \cite{ben}, Bensoussan proved uniqueness of the eigenvector (as
weak solution of the ergodic HJB equation) under assumptions, which translated
in finite dimension imply irreducibility of stochastic matrices
 $P\in\partial f(v)$. Inspired by the results of the present paper, the
first author, Sulem and Taksar~\cite{AST} 
proved uniqueness of the  viscosity
solution of a special ergodic HJB equation.
This yields an example of concrete situation where some non optimal
stationary strategies have several final classes, whereas the optimal ones
have only one final class (translated to our setting, this means that for
some $x\in\R^n$ and $P\in\partial f(x)$, $P$ may have several final classes,
whereas there exists an eigenvector $v$ such that
all elements of $\partial f(v)$ have one final class).

\subsection{From Uniqueness Results to Existence Results}
The uniqueness result that follows from 
Theorem~\ref{t-cst} can be thought of as a partial extension of the
condition of Nussbaum~\cite[Theorem 2.5]{nuss88}: specialized 
to convex monotone 
homogeneous maps $f:\R^n\to\R^n$, the result of Nussbaum shows 
that if $f$ is ${\mathcal C}^1$ and for  all eigenvectors $v$, $f'(v)$ has
only one final class (in this case of course $f$ has a unique critical class),
then the eigenvector of $f$ is unique. The
idea of all these results 
is that the dimension of $\sE(f)$ can be bounded
by looking at ``linearizations'' of $f$ near an eigenvector.

It is instructive to note that the uniqueness
of eigenvectors in $\R^n$ is governed by the same
kind of graph properties as the existence
of eigenvectors, albeit the graphs are different.
For instance, a result of the second author
and Gunawardena~\cite[Theorem 2]{GAUGUN2}
guarantees the existence of an eigenvector for a monotone homogeneous
map which has a strongly connected graph.
Here, the  \new{graph} $\sG(f)$ of a monotone homogeneous
map $f:\R^n\to \R^n$ has nodes $\{1,\ldots,n\}$
and an arc $i\to j$ if $\lim_{\nu \to \infty} f_i(\nu e_j)=+\infty$,
where $e_j$ denotes the $j$-th vector of the canonical basis of $\R^n$.
Another way to guarantee the existence of an eigenvector
is to use the convex
spectral theorem itself, thanks to the following observation
taken from~\cite{GAUGUN2}.
We denote by $\widehat{f}(x)=\lim_{\mu\to\infty}  \mu^{-1}f(\mu x)$
the \new{recession} function of $f$ ($\widehat{f}$ need
not exist when $f$ is monotone and homogeneous,
but it does exist when $f$ is convex).
We have $\widehat f(0)=0$, 
and when $f$ is (additively) homogeneous,
so does $\widehat f$, so that
all points on the diagonal are (trivial) fixed
points of $f$. It is proved in~\cite{GAUGUN2}
that if the recession function of a monotone
homogeneous map $f$ exists and has only fixed
points on the diagonal, then, $f$
has an eigenvector.
Combining this observation with the convex spectral theorem,
we obtain:
\begin{cor}\label{cor-00}
A monotone homogeneous map has an eigenvector
if its recession function exists, is convex,
and has only one critical class.
\end{cor}
If $f$ is a convex monotone homogeneous map,
it is not difficult to see that
the recession function $\widehat{f}$ is exactly the support function
of the domain of the Fenchel transform $f^*$ of $f$ 
(defined  in Section~\ref{cs1} below),
that is $\widehat{f}(x)=\sup_{P\in\dom f^*} Px$.
In this formula, one can replace $\dom f^*$ by its closure $\cl(\dom f^*)$,
which is equal to $\partial \widehat{f}(0)$.
The graph $\sG(f)$ is the union of the graphs of $P\in\dom f^*$,
or equivalently the union of the graphs of $P\in\cl(\dom f^*)$,
whereas the critical graph $\gc(\widehat{f})$
is the union of the final graphs of $P\in\cl(\dom f^*)$.
If $\sG(f)$ is strongly connected, one can
see that $\gc(\widehat{f})$ is also
strongly connected, so that
in the special case of convex monotone homogeneous maps, 
Corollary~\ref{cor-00} is stronger 
than Theorem~2 of \cite{GAUGUN2}
(which however holds in a more general
context).
\subsection{Extension to Convex Monotone Subhomogeneous Maps}
\label{subsec-subhom}
Finally, let us mention two immediate extensions
of the convex spectral theorem.
First, since the map $f\mapsto (x\mapsto -f(-x))$
sends convex monotone homogeneous
maps to concave monotone homogeneous maps,
there is of course a dual 
concave spectral theorem.
Another, more interesting, extension,
is obtained by considering
\new{subhomogeneous}
maps $f$, which satisfy
$f(\lambda +x)\leq \lambda +f(x)$,
for all $\lambda\geq 0$ and $x\in \R^n$.
It is easy to see that a 
monotone map is subhomogeneous if, and
only if, it is nonexpansive for the sup-norm.
To a monotone subhomogeneous map $f:\R^n\to \R^n$,
we associate canonically a monotone homogeneous
map $g:\R^{n+1}\to \R^{n+1}$, \[
g(x,y)=\begin{pmatrix}
y+f(-y+x)\\
y 
\end{pmatrix} \enspace,\quad \forall x\in \R^n,y\in \R 
\]
(this is a nonlinear
extension of the classical way of passing
from a substochastic to a stochastic matrix,
by adding a cemetery state, in this non-linear
context, this construction is due
to Gunawardena and Keane~\cite{gunawardenakeane}).
A vector $z\in \R^n$ is a fixed point of $f$, 
if, and only if $(z,0)$ is an eigenvector
of $g$ (and the eigenvalue is $0$). 
Using this construction, one
translates readily the Convex Spectral Theorem
to a theorem describing fixed point
sets and the asymptotics of the iterates of
convex monotone subhomogeneous maps.
(We might also use this construction, with $-\lambda+f$
instead of $f$, to describe the eigenspace
of $f$ for an additive eigenvalue $\lambda$,
but when $f$ is only monotone and subhomogeneous,
$\lambda$ need not be unique, and it tells little
about the asymptotics of $f^k$, in general.)
For a convex monotone subhomogeneous map $f$ with fixed
point $v$, the \new{critical graph} $\gc(f)$ of $f$ is defined as the union
of the graphs of the matrices $P_{CC}$, where
$P\in \partial{f}(v)$, $C$ is a final class of $P$,
and the $C\times C$ submatrix of $P$, $P_{CC}$, 
is stochastic (when $f$ is homogeneous, 
this property is automatically satisfied). 
Equivalently,
$\gc(f)$ (which can be empty)
is the restriction of $\gc(g)$ to $\{1,\ldots,n\}$.
The notions of critical
classes and cyclicity are defined from $\gc(f)$ as
above.
When $\gc(f)$ is empty,
we have $\nc(f)=\varnothing$, and we take
the convention $\R^\varnothing=\{0\}$, and $\sfc(f)=1$.
\begin{cor}\label{cor-sh}
Let $f:\R^n\to\R^n$ denote a convex monotone
subhomogeneous
map that has a fixed point. Then, all the conclusions
of the Convex Spectral Theorem apply to $f$
and $\lambda=0$. In particular, if $f$ has no critical classes,
then its fixed point is unique.
\end{cor}
{\em Acknowledgements.}\/
The authors thank J. Gunawardena, J.P. Quadrat,
and C. Sparrow, for many useful discussions.

\section{Class Structure of Convex Monotone Homogeneous Maps}\label{cs}
\subsection{Subdifferentials of Convex Monotone Homogeneous Maps}\label{cs1}
We shall first consider {\em scalar}
monotone homogeneous maps $g:\R^n\to\R$
(which satisfy
$x\leq y \implies g(x)\leq g(y)$
for all $x,y\in \R^n$,
and $g(\lambda+x)=\lambda+g(x)$, for all $\lambda\in \R$
and $x\in \R^n$).
The \new{Fenchel transform} of $g$ is the map
$g^*:\R^n\to \Rmin$,
$g^*(p)=\sup_{x\in \R^n}(p\cdot x- g(x))$.
We denote by $\dom g^*=\set{p\in \R^n}{g^*(p)<\infty}$
the \new{domain} of $g^*$, and
by $\proba n=\set{p\in \R^n}{\sum_{1\leq i\leq n}p_i=1,\;
p_1,\ldots,p_n\geq 0}$ the set of \new{stochastic vectors}.
\begin{prop}\label{prop-0}
If $g:\R^n\to \R$ is monotone and homogeneous,
then, $\dom g^*$ is included in $\proba n$.
\end{prop}
\begin{proof}Let $1_n$ denote the vector of $\R^n$ whose entries
are all equal to $1$.
If $g$ is homogeneous, we have for all $p\in \R^n$,
$g^*(p)\geq \sup_{\lambda \in \R} (p\cdot \lambda 1_n-g(\lambda 1_n))
= \sup_{\lambda\in \R} 
(\lambda \sum_{1\leq i\leq n} p_i - \lambda - g(0))
= 
\sup_{\lambda\in \R} 
\lambda (\sum_{1\leq i\leq n} p_i - 1) - g(0)$.
Hence, $g^*(p)<\infty$ implies that $\sum_{1\leq i\leq n} p_i =1$. 
Similarly, if $g$ is monotone, we have $g^*(p)
\geq \sup_{x\leq 0} (p\cdot x - g(x)) \geq 
\sup_{x\leq 0}( p\cdot x- g(0))$.
Hence, $g^*(p)<\infty$ implies that $p_1,\ldots, p_n\geq 0$. 
\end{proof}
\begin{cor}\label{prop-1}
If $g:\R^n\to \R$ 
is monotone, homogeneous,
and convex, then, for all $v\in \R^n$,
$g(v)=\sup_{p\in\dom g^*} (p\cdot v-g^*(p))$, and
the subdifferential 
\begin{eqnarray}
\partial g(v) &\bydef &
 \set{p\in \R^n}{g(x)-g(v)\geq p\cdot (x-v),\;\forall x\in \R^n}\nonumber\\
&= &  \set{p\in \dom g^*}{p\cdot v - g^*(p)=g(v)}\label{e-diesis}
\end{eqnarray}
is a non-empty compact convex subset of $\proba n$.
\end{cor}
\begin{proof}
Since $g$ is convex and takes only finite values, the subdifferential
$\partial g(v)$ is non-empty for any $v\in \R^n$,
and a subdifferential is trivially
closed and convex. Moreover, $\partial g(v)\subset
\dom g^*$ by definition of $\partial g(v)$,
hence, by Proposition~\ref{prop-0}, $\partial g(v)\subset \proba n$.
Finally, $\partial g(v)$ is compact as a closed subset
of the compact set $\proba n$.
\end{proof}
Let $\proba{nn}$ denote the set of $n\times n$ (row-)stochastic matrices.
For any $n\times n$ matrix $P$, we denote by $P_1,\ldots P_n$ the rows of
$P$ and identify $P$ to $(P_1,\ldots , P_n)$, which
amounts to identifying $\R^{n\times n}$
to $(\R^{1\times n})^n$ and $\proba{nn}$ to $(\proba{n})^n$.
If $f:\R^n\to \R^n$ is convex, we set, for all $v\in \R^n$,
\begin{equation}
\partial{f}(v)\bydef \set{P\in \R^{n\times n}}{f(x) - f(v)  \geq P(x-v),\quad \forall x\in \R^n}
\enspace .
\end{equation} 
If $f=(f_1,\ldots,f_n)$, the Fenchel transform of $f$ is the map
$f^*:P\in\R^{n\times n}\mapsto (f_1^*(P_1),\ldots , f_n^*(P_n))\in(\Rmin)^n$,
and its domain is given by 
\[\dom f^*=\set{P\in \R^{n\times n}}{f^*(P)\in\R^n}
=\dom f_1^*\times\cdots\times \dom f_n^*\enspace .\]
Of course 
\begin{equation} 
\partial{f}(v)=\partial{f_1}(v)\times \cdots \times
\partial{f_n}(v) \subset \dom f^* \subset\proba{nn}\label{rect} \enspace .
\end{equation}
when $f$ is monotone and homogeneous.
\subsection{Convex Rectangular Sets of Stochastic Matrices}\label{cs2}
We say that a subset $\sP$ of $\proba{nn}$ is \new{rectangular}
if $\sP=\sP_1\times\cdots\times\sP_n$, for some subsets $\sP_1,\ldots,
\sP_n$ of $\proba{n}$. If $f:\R^n\to\R^n$
 is a convex monotone homogeneous map,
the domain of $f^*$ is rectangular,  and the subdifferential 
$\partial f(v)$ of $f$ at any $v\in\R^n$ is rectangular. 
In this section, we extend the classical definition of graphs
of stochastic matrices
to the case of convex rectangular sets of stochastic matrices. 

Let us first recall some classical notions.
For all graphs $\sG$, and nodes $i,j$ of $\sG$, 
we say that $i$
\new{has access} to $j$, and we
write $i\access j$, if there
is a directed path from $i$ to $j$
in $\sG$, or if $i=j$. We call \new{classes}
of $\sG$ the equivalence classes
for the equivalence relation 
``$(i\access j)\;\text{\rm and}\;(j\access i)$''.
We call \new{strongly connected component} of $\sG$, the restriction
$\sG|_C$ of $\sG$ to a class $C$ of $\sG$, that is the subgraph of $\sG$
with set of nodes $C$ and arcs $i\to j$ when $(i\to j)\in\sG$ and $i,j\in C$.
A graph is strongly connected if it has only one class.
A class $C$ of $\sG$ is \new{final} if
no nodes of $C$ have access to a node
of the complement of $C$.

To any $n\times n$ nonnegative matrix
$P$, we associate the (directed) graph $\sG(P)$
with nodes $1,\ldots,n$ and arcs $i\to j$ when $P_{ij}\neq 0$. 
A matrix is \new{irreducible} if  $\sG(P)$ is strongly connected.
Associating $\sG(P)$ to $P$ allows us to use the graph vocabulary
for nonnegative matrices, for instance the classes of $P$ are by definition
the classes of $\sG(P)$.
If $P$ is a stochastic matrix, a class $C$ of $P$ is final if, and only 
if, the $C\times C$ submatrix of $P$ is stochastic. 
For any stochastic matrix $P$, we denote $\nf(P)$ the union of final classes
of $P$, $\cf(P)$ the set of final classes of $P$ and $\gf(P)$
the subgraph of $\sG(P)$ composed of the nodes and arcs of final classes
of $P$, that is:
\begin{equation}
\gf(P)=\bigcup_{\;\;\;\;F \;\text{\rm final class of}\; P} \sG(P_{FF})
\enspace.
\label{eq-def-gc}
\end{equation}
Here
and below, for all $n\times n$ matrices $P$ and subsets $I$ and $J$
of $\{1,\ldots,n\}$, we denote by $P_{IJ}$
the $I\times J$ submatrix of $P$.
If $\sP$ is a rectangular set of stochastic matrices, we set
\[ \nf(\sP)=\bigcup_{P\in\sP} \nf(P),\quad
\cf(\sP)=\bigcup_{P\in\sP} \cf(P),\quad
\gf(\sP)=\bigcup_{P\in\sP} \gf(P)\, .\]
Hence, $\gf(\sP)$ has $\nf(\sP)$ as set of nodes, and $i\to j$ is an
arc if there exists $P\in\sP$ such that $P_{ij}\not = 0$
and $i$ and $j$ belong to the same final class of $P$.
We say that a class $C\in \cf(P)$ is \new{maximal}
if it is maximal for inclusion.
We have:
\begin{prop}\label{tech}
Let $\sP$ be a convex rectangular set of stochastic matrices.
\begin{enumerate}
\item\label{a1} If $C,C'\in \cf(\sP)$ and $C\cap C'\neq\varnothing$,
then, $C\cup C'\in \cf(\sP)$;
\item\label{a2} The maximal elements of $\cf(\sP)$
are disjoint;
\item\label{a3} There is a matrix $P\in \sP$ such that $\gf(P)=\gf(\sP)$.
The final classes of $P$ are precisely the classes of $\gf(\sP)$.
\item\label{a4} The classes of $\gf(\sP)$ are exactly the maximal
 elements of $\cf(\sP)$.
\end{enumerate}
\end{prop}
\begin{proof}
By assumption,  $\sP=\sP_1\times\cdots\times \sP_n$, where
$\sP_1,\ldots,\sP_n$ are convex subsets of $\proba{n}$.

If $C,C'\in \cf(\sP)$, then, we can find 
 $P,P'\in\sP$ such that $C$ and $C'$ are final classes of $P$
and $P'$, respectively. 
Consider the matrix $P''$
defined by $P''_i= (P_i+P'_i)/2$ if $i\in C\cap C'$,
$P''_i=P_i$ if $i\in C\setminus C'$, 
$P''_i=P'_i$ if $i\in C'\setminus C$,
and $P''_i=P_i$ or $P'_i$ for $i\in \{1,\ldots,n\}
\setminus (C\cup C')$ (the choice has no importance).
Since the $\sP_i$ are convex, $P''_i\in \sP_i$
for all $i$, and by rectangularity, $P''\in\sP$.
To prove the first assertion, it remains to check that 
$C\cup C'$ is a final class of $P''$.
By construction, the 
$(C\cup C')\times (C\cup C')$ submatrix of
$P''$ has row-sum $1$, hence,
we have only to check that all the nodes
of $C\cup C'$ are mutually accessible
in $\sG(P'')$. By definition
of $P''$, any path of $\sG(P)$ (resp. $\sG(P')$) that remains in
$C$ (resp. $C'$) is a path of $\sG(P'')$. Picking
any node $j\in C\cap C'$, we see that
there is a path from any $i\in C\cup C'$ to $j$,
and a path from $j$ to $i$, in $\sG(P'')$, which
shows the first assertion.

The second assertion follows readily from the first.

For the third point, using the method of
the first part of the proof, we shall construct by convexification
a matrix $Q\in \sP$ such that $\gf(\sP)=\gf(Q)$.
For each arc $i\to j$ of $\gf(\sP)$,
there is a matrix $P^{ij}\in \sP$ such that
$i$ and $j$ belong to the same final class $F^{ij}$ of $P^{ij}$,
and $(P^{ij})_{ij}\neq 0$.
For each $k\in \nf(\sP)$,
we consider $\sQ_k=\set{(P^{ij})_{k}}{(i\to j)\in \gf(\sP), k\in F^{ij}}$, 
and we build the stochastic matrix $Q$ whose $k$-th row
is given by $Q_k=|\sQ_k|^{-1}(\sum_{P\in \sQ_k} P)$
for all $k\in \nf(\sP)$
(we denote by $|\cdot|$ the cardinality
of a set), and $Q_k=$ any element of $\sP_k$
if $k\not \in \nf(\sP)$. By convexity and rectangularity of $\sP$,
$Q\in \sP$, hence  $\gf(Q)\subset\gf(\sP)$. 
By construction, each row of $Q$ with index in
$\nf(\sP)$ has sum $1$,
and $\sG(Q)$ contains $\gf(\sP)$.
Moreover, for any arc $(k\to l)$ of $\sG(Q)$ starting at $k\in \nf(\sP)$,
there exists $(i\to j)\in \gf(\sP)$ such that $(P^{ij})_{kl}>0$ and
$k\in F^{ij}$, hence $k$ and $l$ are in $F^{ij}$, which implies 
that $(k\to l)\in \gf(\sP)$.
In particular there are no arcs in $\sG(Q)$ going out from $\nf(\sP)$, and
the restriction of  $\sG(Q)$ to $\nf(\sP)$ is equal to $\gf(\sP)$.
All these remarks imply that $\gf(Q)$ contains  $\gf(\sP)$, so
$\gf(Q)=\gf(\sP)$.

Finally,  any element $F$ of $\cf(\sP)$ (a fortiori,
any maximal element) is such that two nodes of $F$ are connected
by a directed path in  $\gf(\sP)$, so $F$ is included
in a class of $\gf(\sP)$.
Conversely, let $C$ be a class of $\gf(\sP)$, and take 
the matrix $Q\in\sP$ of the third point.
Then $C$ is a final class of $Q$, hence $C\in\cf(\sP)$, so 
$C$ is included in a maximal element of $\cf(\sP)$.
Since the classes of $\gf(\sP)$ are disjoint, we obtain the last point
of the proposition.
\end{proof}

\begin{example}\label{exa1}
Let $\vex(X)$ denote the convex hull of a set $X$, and
consider the convex rectangular set
 $\sP=\sP_1\times\sP_2\times\sP_3\subset\proba{33}$, with
\[
\begin{array}{c}
\sP_1=\vex\{(1,0,0),(1/2,1/2,0),(1/2,0,1/2)\}\enspace ,\\
\sP_2=\vex\{(0,1,0),(2/3,1/3,0)\}\enspace ,\\
\sP_3=\{(0,0,1)\}\enspace .
\end{array}
\]
Since the identity matrix $I$ belongs to $\sP$,
$\{1\},\{2\},\{3\}$ all are elements of $\cf(\sP)$.
Moreover, since
\[
\left(\begin{array}{ccc}
1/2 & 1/2& 0 \\
2/3 & 1/3& 0\\
0 & 0 & 1
\end{array}\right)\in \sP\enspace,
\] 
$\{1,2\}$ is also an element  of $\cf(\sP)$.
It is easy to see that these are the only
elements of $\cf(\sP)$,
so that its maximal elements for inclusion are  
$\{1,2\}$ and $\{3\}$.
The final graph  $\gf(\sP)$ of $\sP$ is the following
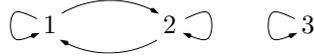
\begin{center}
\begin{tabular}[c]{c}
\input fig1
\end{tabular}
\end{center}
and its classes are $\{1,2\}$ and $\{3\}$.
We shall give, at the end of Section~\ref{piecewise}, an algorithm to determine
the final graph of a convex rectangular sets $\sP$ with a finite number
of extremal points.
\end{example}

\subsection{Critical Graph of Convex Monotone Homogeneous Maps}
\label{ss-crit}
In the sequel, $f$ denotes a convex
monotone homogeneous map that has an eigenvector.
The associated eigenvalue will be denoted by $\lambda$.
As pointed out in the introduction, we may assume $\lambda=0$,
 so that statements and 
proofs will often be given in this case.
We shall prove:
\begin{prop}\label{indep}
The sets  $\nf(\partial f(v))$, $\cf(\partial f(v))$ and the graph 
 $\gf(\partial f(v))$ all are independent of the choice
of the eigenvector $v$ of $f$. 
\end{prop}
Thus, we shall simply write
$\nc(f)$, $\cc(f)$ and $\gc(f)$.
We call \new{critical nodes} of $f$ the elements of $\nc(f)$,
\new{critical graph} of $f$ the graph $\gc(f)$, and
\new{critical classes} of $f$  the classes of  $\gc(f)$.
Combining Proposition~\ref{indep} and Proposition~\ref{tech}, we get:
\begin{cor}\label{crit}
\begin{enumerate}
\item\label{b2} For any eigenvector $v$ of $f$,
there is a matrix $P\in \partial f(v)$ such that $\gc(f)=\gf(P)$.
In particular, the final classes of $P$ are precisely the  critical classes of $f$.
\item\label{b1} The critical classes of $f$ are exactly 
the maximal elements of $\cc(f)$.
\end{enumerate}
\end{cor}

The graph $\gc(f)$ generalizes the classical critical graph that appears
in max-plus algebra. Let us recall that if $A$ is a $n\times n$ matrix
with entries in $\Rmax$, the \new{critical graph} $\gc(A)$ is the union of the
circuits $(i_1,\ldots,i_k)$ that realize the maximum in the formula:
\[ \rho(A)= \max_{1\leq k\leq n}\;  \max_{i_1,\ldots, i_k}
\frac{A_{i_1i_2}+\cdots+A_{i_ki_1}}{k}\enspace . \]
If $A$ has at least one finite entry per row, the max-plus linear map $f_A$
with matrix $A$, $(f_A(x))_i=\max_{1\leq j\leq n} (A_{ij}+x_j)$,
sends $\R^n$ into $\R^n$. Of course, $f_A$ is monotone, homogeneous and
convex.
If $v\in\R^n$ is a max-plus eigenvector of $A$, 
with associated eigenvalue $\lambda$ (i.e. if $f_A(v)=\lambda+v$),
we define the \new{saturation graph}
of $A$ and $v$, $\sat(A,v)$, as the union of arcs $i\to j$ such that
$\lambda+ v_i= A_{ij}+v_j$.
\begin{prop}\label{max-crit}
If $A\in (\Rmax)^{n\times n}$ has at least one finite entry per row, and if
$v\in\R^n$ is an eigenvector of $A$, then $\gc(f_A)=\gc(A)$ coincides with the
union of strongly connected components of $\sat(A,v)$.
\end{prop}
\begin{proof}
The fact that $\gc(A)$ coincides with the
union of strongly connected components of $\sat(A,v)$ is a consequence of 
the max-plus spectral theorem (see~\cite[\S3.2.4 and 3.7]{BCOQ}).
The rest of the assertion follows from the identity
$ (\partial f_A(v))_i=\vex\set{\delta_j}{(i,j)\in \sat(A,v)}$,
where $\delta_j\in\proba n$ denotes the Dirac probability measure at $j$.
Indeed, from this, one deduces that, for any $P\in\partial f_A(v)$,
$\sG(P)\subset \sat(A,v)$, so that $\gf(P)$ is included in the union of
strongly connected components of $\sat(A,v)$, that is $\gc(A)$.
Hence, $\gc(f_A)\subset \gc(A)$.
Conversely, for any circuit $c$ in $\gc(A)\subset \sat(A,v)$,
one can construct $P\in\partial f_A(v)$ such that $P_{ij}=1$ for any
arc $i\to j$ of $c$. The circuit $c$ , which is a final class of $P$,
is included in $\gc(f_A)$. Thus, $\gc(A)\subset \gc(f_A)$.
\end{proof}
(We shall extend this property to the stochastic control context in
Proposition~\ref{t-o}.)

Proposition~\ref{indep} relies on classical
facts of Perron-Frobenius theory. We first need the characterization
of supports of invariant measures.
Recall that an \new{invariant measure}
of a stochastic matrix $P$ is a stochastic
(row) vector $m$ such that $mP=m$.
The set of non-zero entries of $m$
is the \new{support} of $m$.
The following result is standard:
\begin{prop}\label{propf}
Let $P$ be a $n\times n$ stochastic matrix.
A subset $F$ of $\{1,\ldots,n\}$ is a union
of final classes of $P$ if, and only if, $P$
has an invariant measure $m$ with support
$F$.\qed
\end{prop}
The main ingredient of the proof of Proposition~\ref{indep}
is a classical property of super-harmonic vectors
(a vector $v$ is \new{super-harmonic} 
for a stochastic matrix $P$
if $Pv\leq v$), which is itself a consequence of the Perron-Frobenius theorem.
\begin{lemma}[Discrete Maximum Principle]
\label{disc}
If $P$ is a $n\times n$ stochastic matrix,
and if $Pv \leq v$, then, the restriction of $v$ 
to any final class of $P$ is a constant vector,
and $Pv=v$ holds on the union $\nf(P)$ of the final classes of $P$.
Moreover, the minimum of $v$ is attained on a final class.
\end{lemma}
\begin{proof}
Let $C$ denote a final class of $P$. It follows from Proposition~\ref{propf}
that there exists an invariant measure $m$ of $P$ with support $C$.
We have, $0 \leq m(v-Pv)= mv-mv=0$,
and since $0=m(v-Pv)=\sum_{i\in C}m_i(v_i-(Pv)_i)$
is a sum of nonnegative terms, we must have $v_i=(Pv)_i$
for all $i\in C$. This means that the restriction of $v$
to $C$ is a right-eigenvector of the $C\times C$ submatrix
of $P$. Since this submatrix is stochastic
and irreducible, this eigenvector
must be constant, which shows the first part
of the lemma. 

To show the second part,
we introduce $F=\nf(P)$ and $T=\{1,\ldots,n\}\setminus F$. We denote by $Q$ the
$T\times T$ submatrix of $P$, and by $w$ the restriction
of $v$ to $T$. Let $\mu$ denote the minimum of the entries
of $v$ restrained to $F$. Since $P1_n=1_n$, $v-\mu \geq P(v-\mu)$,
which allows us to assume that $\mu=0$.
Then, $v\geq Pv$ and $v\geq 0$ on $F$ yield $w\geq Qw$, and since
$Q$ has nonnegative entries,
$w\geq Qw\geq Q^2w\geq \ldots \geq Q^kw$,
for all $k\geq 1$. Since $Q$ has spectral radius
strictly less than $1$,
we get $w\geq \lim_k Q^kw=0$.
\end{proof}

We finally prove Proposition~\protect\ref{indep}.
Let $v,v'$ denote two eigenvectors
of $f$. Let $C\in \cf(\partial f(v))$
and let $P\in \partial f(v)$
for which $C$ is a final class. We get
$v'-v=f(v')-f(v)\geq P(v'-v)$.
It follows from Lemma~\ref{disc} that $(v'-v)_i=P_i(v'-v)$
holds for all $i\in C$.
But $f_i(v')-f_i(v)=(v'-v)_i=P_i(v'-v)$
and $P_i\in \partial{f}_i(v)$ imply that $P_i\in \partial{f}_i(v')$
(indeed, $f_i(x)-f_i(v')=f_i(x)-f_i(v)+ f_i(v)
-f_i(v') \geq P_i(x-v) -P_i(v'-v)
\geq P_i(x-v')$ holds for all $x\in \R^n$).
Now, pick any $Q\in \partial f(v')$, and
consider the matrix $R$ such that $R_i=P_i$
for all $i\in C$ and $R_i=Q_i$ if $i\in \{1,\ldots,n\}\setminus C$.
By rectangularity of $\partial{f}(v')$, $R\in \partial{f}(v')$. Since
$C$ is a final class of $R$, $C\in \cf(\partial{f}(v'))$.
This shows that  $\cf(\partial{f}(v))\subset \cf(\partial{f}(v'))$.
Since $R_{CC}$ and $P_{CC}$ coincide, this also shows that
 $\gf(\partial{f}(v))\subset \gf(\partial{f}(v'))$, and thus
 $\nf(\partial{f}(v))\subset \nf(\partial{f}(v'))$.
\qed

\begin{example}\label{exa}
Let $f:\R^3\to \R^3$,
\begin{equation}
f(x)=\left(
\begin{array}{c}
x_1 \vee \frac{1}{2}(x_1+x_2) \vee (-3+x_3) \vee \frac{1}{2} (x_1+x_3)\\
x_2 \vee \frac{1}{3}(2x_1+x_2)\vee (-1+x_1)\\
x_3 \vee (-2+ x_1)
\end{array}
\right)\enspace, 
\label{eq-f1}
\end{equation}
which is such that $f(0)=0$.
We have $\partial{f}(0)=\sP$, with $\sP$ defined in Example~\ref{exa1}
(that is $\partial{f}_i(0)=\sP_i$ for $i=1,2,3$).
So, the critical graph $\gc(f)$ of $f$ is exactly the final graph of 
$\sP$ shown in Example~\ref{exa1}, and the critical classes of $f$ 
are $\{1,2\}$ and $\{3\}$.
Let us illustrate the fact that $\gf(\partial f(v))$ is independent
of the choice of the eigenvector $v$ by considering
the eigenvector $v=(0,0,-2)$. By comparison with $\partial{f}(0)$,
only the first and third rows of $\partial{f}(v)$ are changed,
namely $\partial{f}_1(v)=\vex\{(1,0,0),(1/2,1/2,0)\}$ and
$\partial{f}_3(v)=\vex\{(0,0,1),(1,0,0)\}$.
Although $\partial{f}(v)\neq \partial{f}(0)$,
the final classes
of the matrices $P\in \partial{f}(v)$ are
the same as the ones of the matrices $P\in \partial{f}(0)=\sP$.
\end{example}

\section{Structure of Eigenspaces}
\label{eigen}
In order to make more apparent the proof
idea, we first show a simple result.
\begin{theorem}\label{theo-1}
The eigenvector of a convex monotone homogeneous
map with a unique  critical class
is unique, up to an additive constant.
\end{theorem}
\begin{proof}
Let $C$ denote the  critical class of $f$,
and let $v,v'$ be two eigenvectors of $f$.
Using Assertion~\ref{b2} of Corollary~\ref{crit},
we get matrices $P\in \partial f(v)$
and $P'\in \partial f(v')$ such that $C$
is the unique final class of $P$ and $P'$.
Since $P\in \partial f(v)$,
$v'-v=f(v')-f(v)\geq P(v'-v)$, hence,
by Lemma~\ref{disc}, $v'-v$ is
constant on $C$, and it attains its minimum
in $C$. Exchanging the roles of $v'$ and $v$,
we see that $v'-v$ attains its maximum in $C$.
Therefore, $v'-v$ is a constant vector.
\end{proof}

\begin{example}
Let us use Theorem~\ref{theo-1} to show that
the map $f:\R^3\to \R^3$,
\[
f(x)=\left(
\begin{array}{c}
\frac{1}{2}(x_1+x_2)\\
\log(e^{x_2} + 8 e^{x_3})\\
x_1 \vee (-1+x_3)
\end{array}
\right)
\]
has a unique eigenvector, 
up to an additive constant.
Since the graph $\sG(f)$ of $f$, which is equal to 
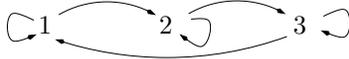
\begin{center}
\begin{tabular}[c]{c}
\input fig2 
\end{tabular}\enspace
\end{center}
is strongly connected, we known, thanks to \cite[Theorem 2]{GAUGUN2},
that $f$ has an eigenvector. Indeed, $f(v)=\lambda + v$, where
$\lambda=\log 2$, and $v=(\log 2,3\log 2,0)$. 
We have
\[
\partial f(v)=\{f'(v)\}\,\;
\text{\rm where}\;
f'(v)= \left(
\begin{array}{ccc}
1/2 & 1/2 & 0\\
0 &1/2& 1/2\\
1 & 0& 0 
\end{array}
\right)\enspace .
\]
Thus, the critical graph $\gc(f)=\gf(f'(v))=\sG(f'(v))$ is 
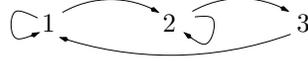
\begin{center}
\begin{tabular}[c]{c}
\input fig3
\end{tabular}
\end{center}
(only the arc $3\to 3$ vanished, by comparison
with $\sG(f)$). Since $\gc(f)$ is strongly connected,
$f$ has a unique critical
class (namely $\{1,2,3\}$), and, by Theorem~\ref{theo-1},
the eigenvector of $f$ is
unique (up to an additive constant).
Equivalently, the map $F:(\Rp)^3\to (\Rp)^3$, $F=\exp \comp f\comp \log$,
\[
F(x)=\left(
\begin{array}{c}
\sqrt{x_1x_2}\\
x_2 + 8 x_3 \\
x_1 \vee (x_3/e) 
\end{array}
\right) ,
\]
has a unique (multiplicative) eigenvector,
up to a constant, namely, $(2,8,1)$.
\end{example}

To extend Theorem~\ref{theo-1} to the case
of several  critical classes, we shall
use the following fundamental observation,
which establishes the existence
of a monotone homogeneous projector $(-\lambda+f)^\omega$ 
from the \new{super-eigenspace} of $f$, 
$\sE^+(f)=\set{x\in \R^n}{f(x)\leq \lambda+x}$,
to the \new{eigenspace} $\sE(f)=\set{x\in \R^n}{f(x)=\lambda+x}$.
\begin{lemma}[Spectral Projector]\label{lem-projection}
If $f:\R^n\to\R^n$ is a convex monotone homogeneous map with eigenvalue $0$,
and if $z\in \R^n$ is such that $f(z)\leq z$, then
\[f^\omega(z)\bydef \lim_{k\to\infty} f^k(z)
\]
is an eigenvector of $f$, which coincides with $z$ on
the union $\nc(f)$ of  critical classes of $f$.
The map $f^\omega:\sE^+(f)\to \sE(f)$ is monotone, homogeneous, convex, onto, 
and satisfies $(f^\omega)^2=f^\omega$.
\end{lemma}
\begin{proof}
By monotonicity of $f$, the sequence
$z_k$ defined by $z_0=z$, $z_k=f(z_{k-1})$ for $k\geq 1$,
is nonincreasing.
Picking any eigenvector $v$ of $f$ and using
the nonexpansiveness of $f$,
we get $|z_k-v|=|f^k(z_0)-f^k(v)|\leq |z_0-v|$.
Therefore, the sequence $z_k$
which is bounded and nonincreasing converges
to a limit $f^\omega(z)$,
which by continuity of $f$ is an eigenvector of $f$.
To show that $f^\omega(z)$ coincides with $z$ on 
$C=\nc(f)$, we observe
that, for all $u\in \R^n$,
\begin{equation}
f(u)\leq u\implies f(u)=u \;\text{\rm on}\; C \enspace.
\label{eq-saturate}
\end{equation}
Indeed, arguing as in the proof
of Theorem~\ref{theo-1},
we can find, by Corollary~\ref{crit}, a matrix $P\in \partial{f}(v)$
which has $C$ as union of final classes.
Then, $u-v\geq f(u)-f(v)\geq P(u-v)$,
and by Lemma~\ref{disc}, $u-v=P(u-v)$ on $C$,
which implies that $u-v=f(u)-f(v)$ on $C$,
and since $v=f(v)$, we get~\eqref{eq-saturate}.
Applying~\eqref{eq-saturate} to
the inequalities $f(z_0)\leq z_0$, $f(z_1)\leq z_1$, \ldots,
we get $z_1=f(z_0)=z_0$ on $C$, $z_2=f(z_1)=z_1$ on $C$, \ldots,
hence $z_k=z_0$ on $C$ for all $k\geq 0$, which implies
that $f^\omega(z)=\lim_k z_k=z$ on $C$.
Finally, $f^\omega:\sE^+(f)\to\sE(f)$ is onto
since $f^\omega(x)=x$ when $x\in \sE(f)$,
and $f^\omega$, which is a pointwise
limit of the monotone, homogeneous, and convex maps $f^k$
is also monotone, homogeneous, and convex.
\end{proof}
The super-eigenspace $\sE^+(f)$ will be useful in the subsequent
proofs, because $\sE^+(f)$ is a ``more regular'' object
than $\sE(f)$. In particular, thanks to the monotonicity of $f$,
$\sE^+(f)$ is an inf-subsemilattice of $(\R^n,\leq)$,
and, when $f$ is convex, $\sE^+(f)$ is convex.
If $C$ is any subset of $\{1,\ldots,n\}$,
we denote by $\R^C$ the set of vectors indexed
by the elements of $C$ and by $\sfr_C$ the restriction map
$\R^n\to\R^C$, $\sfr_C(x)_i=x_i$, for all $i\in C$.
\begin{theorem}\label{thmain}
Let $C_1,\ldots, C_s$ denote the 
critical classes of a convex monotone homogeneous map $f$,
and let $C=C_1\cup \ldots \cup C_s=\nc(f)$.
Then, any two eigenvectors $v,v'$ of $f$ satisfy:
\begin{enumerate}
\item $v-v'$ is constant on each $C_i$;
\item if $v=v'$ on $C$, then $v=v'$.
\end{enumerate}
Moreover,
\begin{enumerate}
\setcounter{enumi}{2}
\item\label{i-2} the restriction $\sfr_C$
sends bijectively $\sE(f)$ to a convex inf-subsemilattice of $\R^C$,
denoted by $\ec(f)$,
and the inverse map $\sfr_C^{-1}:
\ec(f)\to \sE(f)$ is monotone and homogeneous.
\end{enumerate}
\end{theorem}
\begin{proof}
Let $v,v'$ be two eigenvectors of $f$. 
Arguing as in the proof of Theorem~\ref{theo-1},
we get two matrices
$P\in \partial f(v)$ and $P'\in\partial f(v')$
which both have exactly the final classes $C_1,\ldots,C_s$.
Since $v'-v\geq P(v'-v)$, we get from Lemma~\ref{disc}
that $v'-v$ is constant on each final class $C_i$,
which proves the first assertion. Moreover, if $v'-v\geq 0$
on $C$, the vector $v'-v$, which, by Lemma~\ref{disc},
attains its minimum on a final class, must be nonnegative. 
By symmetry, if $v-v'\geq 0$ on $C$, $v-v'$ must be nonnegative.
Thus, $v=v'$ on $C$ implies $v=v'$, which shows the second assertion
of the theorem. This assertion implies
that $\sfr_C$ is injective $\sE(f)\to \ec(f)$.
The inverse map $\sfr_C^{-1}$ is trivially
homogeneous. It is also monotone: if $w'=\sfr_C(v')
\geq w=\sfr_C(v)$, $v'-v\geq 0$ on $C$, and we just
saw that this implies that $v'-v\geq 0$. Thus, $w'\geq w\implies 
\sfr_C^{-1}(w')\geq \sfr_C^{-1}(w)$.

Take now $v,v'$ two eigenvectors of $f$, and
assume that the eigenvalue of $f$ is $0$.
By monotonicity of $f$, $f(v\wedge v')\leq v\wedge v'$,
and by Lemma~\ref{lem-projection},
$f^\omega(v\wedge v')$ is an eigenvector which coincides
with $v\wedge v'$ on $C$. This shows that $\ec(f)$,
equipped with the ordering $\leq$, is an inf-semilattice,
in which the greatest lower bound coincides
with that of $\R^C$.

Finally, let $0<\alpha<1$, and $\alpha'=1-\alpha$. We have,
by convexity of $f$,
$f(\alpha v+ \alpha'v')\leq \alpha f(v)+\alpha'f(v')
=\alpha v+\alpha 'v'$, and arguing as for $v\wedge v'$,
we derive from this that $\ec(f)$ is convex.
\end{proof}

Of course, Theorem~\ref{thmain}
implies that the eigenspace $\sE(f)$ is connected
($\sE(f)$ is the image of the convex set $\ec(f)$
by the continuous map $\sfr_C^{-1}$). More generally,
Bruck~\cite{bruck} proved that the set of fixed
points of a nonexpansive map in a Banach space
is a retract of the whole space, and, a fortiori,
is connected. In the special case of convex
monotone homogeneous maps, we have
more, since Theorem~\ref{thmain} shows that the eigenspace $\sE(f)$
is isomorphic (by monotone homogeneous bijections)
to the convex set $\ec(f)$.
This makes the following definition natural.
\begin{defi}
The \new{dimension} $\dim \sE(f)$ of the eigenspace $\sE(f)$
of a convex monotone homogeneous map $f$ 
is the dimension (as a convex set) of its restriction
 $\ec(f)$ to the union of  critical classes.
\end{defi}
(Recall that the dimension of a convex set is
defined as the dimension of its affine hull).
A perhaps more convenient way to describe the geometry
of $\sE(f)$ is to introduce the notion
of \new{section}, which is any set
$S\subset \nc(f)$
which meets each critical class at exactly one point.
We denote by $\sfm(f)$ the number of critical
classes of $f$.  Then, paraphrasing Theorem~\ref{thmain},
we get:
\begin{cor}\label{cor-summarize}
Let $f$ denote a convex monotone homogeneous map $f$
with section $S$. 
Then, $\ec(f)$ and
$\sfr_S(\sE(f))$ both are convex sets of dimension $\dim \sE(f)\leq \sfm(f)$.
\end{cor}
A natural question is to characterize the dimension
of $\sE(f)$. We shall see that the equality $\dim \sE(f)=\sfm(f)$ holds,
except in ``singular'' cases.
Let us first exhibit a regular case.
\begin{example}
For the map $f$ given in~\eqref{eq-f1}, we saw in Example~\ref{exa}
that the  critical classes of $f$ are $\{1,2\}$
and $\{3\}$. In particular, $\sfm(f)=2$. Since $v=0$ is an eigenvector,
the first statement of Theorem~\ref{thmain} shows
that any eigenvector is constant on $\{1,2\}$.
Then, an elementary computation shows that
$\sE(f)=\set{x\in \R^3}{f(x)=x}=
\set{x\in \R^3}{x_1=x_2\geq x_3 \geq -2+ x_1}$.
Thus $\dim \sE(f)=2=\sfm(f)$.
\end{example}
\begin{example}
It is instructive to apply
Theorem~\ref{thmain} or Theorem~\ref{theo-1}  to the
case of the map 
$f(x)=\log M\exp(x)$, i.e.,
$f_i(x)=\log(\sum_j M_{ij}\exp(x_j))$,
where $M$ is an irreducible nonnegative matrix $M$.
We have $\partial{f}(x)= \{P(x)\}$,
where $P(x)=f'(x)$ is the $n\times n$ matrix with entries
$P(x)_{ij}=(\sum_k M_{ik}\exp(x_k))^{-1}M_{ij}\exp(x_j)$.
Since $M$ is irreducible, $P(x)$ is irreducible,
for all values of $x$. Therefore, $f$ has a unique 
 critical class, and Theorem~\ref{theo-1}
tells that its eigenvector is unique, up to an additive constant,
or, equivalently, that the matrix $M$ has a unique positive
eigenvector, up to a multiplicative constant.
This is a (complicated) way to derive the uniqueness part of 
Perron-Frobenius theorem for irreducible  nonnegative matrices,
from the Perron-Frobenius theorem for irreducible stochastic matrices.
\end{example}
The next example illustrates the difference between
regular and singular cases.
\begin{example}\label{moreau}
Let $h$ denote any convex map $\R\to \R$ whose subgradients
are between $0$ and $1$. The map $f:\R^2\to \R^2$,
\begin{equation}
f(x)= \left(\begin{array}{c}
x_1 + h(x_2-x_1)\\
x_2 + h(x_1-x_2) 
\end{array}\right)\enspace,
\label{eq-generic}
\end{equation}
is monotone, homogeneous, convex, and $f(0)=h(0)+0$, 
which means that $0$ is an eigenvector of $f$ for the eigenvalue $h(0)$.
If ($h(t)=h(-t)\implies t=0$), then $0$ is the unique eigenvector
of $f$, up to an additive constant. This condition
is satisfied by the maps $t\mapsto\log(1+e^t)$
and $t\mapsto 0\vee t$. It is also satisfied by the Moreau-Yoshida
regularization, $h_\epsilon$, of $t\mapsto 0\vee t$,
which, for all $\epsilon>0$, is defined by
$h_\epsilon(t) =0$ for $t\leq 0$,
$h_\epsilon(t)=t^2/2\epsilon$ for $0\leq t\leq \epsilon$, 
and $h_\epsilon(t)=t-\epsilon/2$ for $t\geq \epsilon$.
Let $g$, $k$, and $k^\epsilon$  
denote the maps $f$ obtained by replacing $h$ by 
$t\mapsto \log (1+e^t)$, $t\mapsto 0\vee t$, and
$h_\epsilon$ in~\eqref{eq-generic}, respectively.
We have $\partial{g_1}(0)=\{(g_1)'(0)\}=\{(1/2, 1/2)\}$,
$\partial{k_1}(0)=\vex\{(1,0),(0,1)\}$, 
and $\partial{k^\epsilon_1}(0)=\{(k^\epsilon_1)'(0)\}=\{(1, 0)\}$.
Using the symmetry between the first and the second coordinate
of $f$, we get the following critical graphs
\begin{center}
\begin{tabular}[c]{c}
\input fig4
\end{tabular}\enspace .
\end{center}
Thus, $\dim\sE(g)=\sfm(g)=1$, 
$\dim\sE(k)=\sfm(k)=1$, 
but $\dim \sE(k^\epsilon)=1<\sfm(k^\epsilon)=2$.
This discrepancy should  be intuitively clear
by looking at the following graphs:
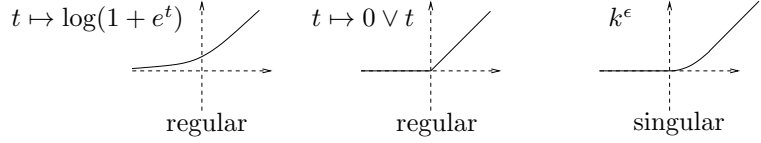
\begin{center}
\begin{tabular}[c]{c}
\input fig5
\end{tabular} \enspace .
\end{center}
In the case of $f=k^\epsilon$,
Theorem~\ref{thmain} fails to characterize
the dimension of $\sE(f)$ because the subdifferential of $f$ 
does not give enough information on the local behavior of $f$
near the eigenspace. In such cases, one needs
to consider terms of higher order
in the local expansion of $f$ to establish
the uniqueness of the eigenvector.
However, we shall see in section~\ref{piecewise}
that the equality $\dim \sE(f)=\sfm(f)$ does
hold when $f$ is piecewise affine.
\end{example}
The following last example shows that $\sE(f)$
need not be convex or an inf-subsemilattice of $\R^n$, even if
its restriction $\ec(f)$ is
a convex inf-subsemilattice of $\R^C$.
\begin{example}\label{ex-f3}
Consider $f:\R^3\to\R^3$,
\[
f(x)= \left(\begin{array}{c}
x_2\vee x_3\\
x_2\\
x_3
\end{array}\right)
\enspace .
\]
We have $\sE(f)=\set{(x_2\vee x_3, x_2,x_3)}{x_2,x_3\in \R}$,
$f$ has two  critical classes $C_1=\{2\}$,
$C_2=\{3\}$, so $\dim \sE(f)=\sfm(f)=2$.
Let $u=(1,1,0)$, $v=(1,0,1)$.
Since $u,v\in \sE(f)$ but
$(u+v)/2=(1,1/2,1/2)\not\in \sE(f)$,
$\sE(f)$ is not convex,
and since $u\wedge v=(1,0,0)\not\in \sE(f)$,
$\sE(f)$ is not an inf-subsemilattice of $\R^C$.
\end{example}
Finally, we note that fixed points sets of monotone
nonexpansive maps $\R^n\to \R^n$ are {\em always} lattices:
the interesting lattice 
statement in assertion~\ref{i-2} of Theorem~\ref{thmain}
is that the inf law of $\sE^c(f)$ {\em coincides} with the
usual inf law in $\R^C$. 
\begin{prop}
Let $f:\R^n\to \R^n$ denote a monotone 
map that is nonexpansive (for any norm) and that has a fixed point.
Then, the fixed point set $\sE(f)=\set{x\in \R^n}{f(x)=x}$,
equipped with the standard order relation, is a lattice, 
in which the sup and inf laws, $\vee_{f}$ and $\wedge_f$, respectively, are
given by 
\begin{equation}
x\vee_f y =\lim_{k\to\infty} f^k(x \vee y),\qquad
x\wedge_f y =\lim_{k\to\infty} f^k(x \wedge y)\enspace. 
\label{e-lattice}
\end{equation}
\end{prop}
\begin{proof}
If $x,y\in \sE(f)$, then by monotonicity of $f$,
$f(x\wedge y)\leq f(x)\wedge f(y)=x\wedge y$.
Arguing as at the beginning of the
proof of Lemma~\ref{lem-projection},
we get that the limit $w=\lim_{k\to\infty} f^k(x\wedge y)$ exists
and satisfies $f(w)=w$,  $w\leq x$, and $w\leq y$.
Moreover, if $z$ is an arbitrary element in $\sE(f)$ such that 
$z\leq x$, and $z\leq y$, we
have $z=f(z)=\cdots = f^k(z)\leq f^k(x\wedge y)$
for all $k$, so that $z\leq w$. This shows that
$w=\lim_{k\to\infty} f^k(x\wedge y)=x\wedge_f y$. The dual
argument shows the first equality in~\eqref{e-lattice}.
\end{proof}

\section{Critical Graph of $f^k$}
\label{secttang}
In this section, we establish the following result which is central in the 
proof of the cyclicity theorem.
If $\sG$ is a graph, we call $\sG^k$ the graph with same nodes as $\sG$
and arcs $(i\to j)$ when there is a directed path $i\to i_1\to\cdots \to i_k=j$
with length $k$ in $\sG$.

\begin{theorem}\label{ncf}
Let $f$ denote a convex monotone homogeneous map $\R^n\to\R^n$
that has an eigenvector. Then, for all $k\geq 1$,
$\gc(f^k)=\gc(f)^k$, in particular $\nc(f^k)=\nc(f)$.
\end{theorem}

The proof of Theorem~\ref{ncf} needs tools and results,
of independent interest, which involve
one-sided directional derivatives, additive recession functions, 
and a nonlinear generalization of invariant measures.
Let us first show the linear version of Theorem~\ref{ncf}.

\begin{prop}\label{cor-c}
If $P$ is a stochastic matrix,
then, for all $k\geq 1$, $\gf(P^k)=\gf(P)^k$, in particular 
$\nf(P^k)=\nf(P)$.
\end{prop}
\begin{proof}
For all nonnegative matrices $P$, we have $\sG(P^k)=\sG(P)^k$.
Since $\gf(P)$  is equal to $\sG(Q)$,  where $Q$ is the restriction of $P$
to $\nf(P)$ and since the restriction of $P^k$ to $\nf(P)$ is 
equal to $Q^k$ (by definition of final classes), it
is enough to prove that $ \nf(P^k)=\nf(P)$.

To show this, we shall use the following consequence of Proposition~\ref{propf}:
for any stochastic matrix $P$,
the support of an invariant measure $m$ of $P$ is a subset of $\nf(P)$,
and there exists an invariant measure $m$ of $P$ with support $\nf(P)$.

Consider first  an invariant measure  $m$ of $P$ with support $\nf(P)$.
The vector  $m$ is also an invariant measure of $P^k$, 
which implies that its support $\nf(P)$ is included in $\nf(P^k)$.
Conversely, consider an invariant measure  $m$ of $P^k$ 
with support $\nf(P^k)$.
We get $m( I+\cdots + P^{k-1}) P= m (P+\cdots + P^k)= m(I+\cdots + P^{k-1})$,
so  $m_k=\frac{1}{k} m( I+\cdots + P^{k-1})$  is an invariant measure of $P$.
Hence, the support of $m_k$ is included in $\nf(P)$.
Since $m_k\geq \frac{1}{k} m\geq 0$, due to the nonnegativeness of $P$ 
and $m$, the support of $m_k$ 
contains that of $m$, which leads to $\nf(P^k)\subset \nf(P)$.
\end{proof}

For any convex map $f:\R^n\to\R$
and any $v\in \R^n$, we denote by $f'_v(y)$
the (one-sided) \new{directional derivative} of $f$ at $v$
with respect to $y$:
\[
f'_v(y)=\lim_{\varepsilon\to 0^+} \frac{f(v+\varepsilon y)-f(v)}{\varepsilon}
\enspace. \]
The map $f'_v$ is well-defined, finite, convex, multiplicatively
homogeneous, and
\begin{equation}\label{propfp}
f'_v(y)=\sup_{p\in \partial f(v)}p\cdot y \enspace,
\end{equation}
see \cite[Theorems~23.1 and 23.4]{ROCK}.
When $f$ is monotone and additively homogeneous,
so does $f'_v$. In the sequel, we shall say that
a map is \new{bihomogeneous} if it is both additively
and multiplicatively homogeneous.

The definition and properties of $f'_v$ can be extended to convex maps
$f:\R^n\to\R^m$, with $f'_v=((f_1)'_v, \ldots ,(f_m)'_v):\R^n\to\R^m$. 
We have the following chain rules:
\begin{lemma}\label{chain}
Let $f:\R^n\to\R^m$ and $g:\R^m\to\R^\ell$ be two convex maps.
Assume that $g$ is monotone. Then, $g\circ f:\R^n\to\R^\ell$ is convex, and,
for all $v\in\R^n$,
\begin{eqnarray}
(g\circ f)'_v&=&g'_{f(v)}\circ f'_v\enspace, 
\label{cr-1}\\
\partial (g\circ f)(v) &=&
\vex \big(\partial g(f(v)) \partial f(v)\big) 
\enspace .
\label{cr-2}
\end{eqnarray}
\end{lemma}
(Recall that $\vex$ denotes the convex hull. In~\eqref{cr-2},
$\partial g(f(v)) \partial f(v)\bydef\set{PQ}{P\in \partial g(f(v)) ,\;Q\in \partial{f}(v)}$.)
\begin{proof}
The convexity of $g\circ f$ is immediate,
and Eqn~\eqref{cr-1} follows easily from
the fact that any finite convex function is locally
Lipschitz continuous. 

To show~\eqref{cr-2}, we first note that if $D$
is a rectangular subset of $\R^m$, i.e. if $D$ is of the form
$D_1\times \cdots \times D_m$, and if $h:\R^m\to \R^p$ is monotone
and satisfies $h(\lim_k x_k)=\lim_k h(x_k)$
for all nondecreasing sequences $x_k$, 
then, 
\begin{equation}
h(\sup_{d\in D} d) =\sup_{d\in D} h(d) \enspace .
\label{e-hsup}
\end{equation}
Indeed, the inequality $\geq$ in~\eqref{e-hsup}
follows from the monotonicity of $h$, and, thanks
to the rectangularity of $D$, we can find a nondecreasing
sequence $d_k\in D$ such that $\sup_{d\in D} d= \lim_{k} d_k$,
hence, $h(\sup_{d\in D} d)= h(\lim_k d_k)=\lim_k h(d_k)
\leq \sup_{d \in D} h(d)$, which shows~\eqref{e-hsup}.

Applying~\eqref{e-hsup} to $D=\partial{f(v)}y$
and $h(x)=Px$, with $P\in \partial{g}(f(v))$,
we get from~\eqref{cr-1} and~\eqref{propfp}, 
\begin{eqnarray}
(g\circ f)'_v(y)=g'_{f(v)}\circ f'_v(y)
&=& \sup_{P\in \partial g(f(v))} P (\sup_{Q\in \partial{f}(v)} Q y)\nonumber\\
&=& \sup_{R\in \partial g(f(v)) \partial{f}(v)} R y  \enspace .\label{cr-3}
\end{eqnarray}
Let us denote by $\delta_X^*:\R^n\to \R\cup\{+\infty\}$
the \new{support function}
of a subset $X\subset \R^n$, which is defined
by $\delta_X^*(p)=\sup_{x\in X} 
p\cdot x$. We first assume that $\ell=1$.
Then,~\eqref{cr-3} and property~\eqref{propfp}
show that $(g\circ f)'_v =\delta^*_{\partial(g\circ f)(v)}=
\delta^*_{\partial{g}(f(v))\partial f(v)}$.
By Legendre-Fenchel duality,
two subsets of $\R^n$
have the same support function, if, and only if, they
have the same closed convex hull
(see~\cite[Corollary~13.1.1]{ROCK}). Since
$\partial(g\circ f)(v)$ is closed and convex,
$\partial(g\circ f)(v) = \vexb \left(\partial{g}(f(v))\partial f(v)\right)$
where $\vexb$ denotes the closed convex hull of a set.
But $\partial{g}(f(v))\partial f(v)$, which is
the image by the continuous map $(P,Q)\mapsto PQ$ of 
the product of two compact sets, is compact, 
and the  convex hull of a compact
subset of $\R^n$, which is compact, coincides
with its closed convex hull 
(see for instance~\cite[Corollary~5.18]{aliprantis}), 
therefore
$\vexb  \left(\partial{g}(f(v))\partial f(v)\right)=
\vex  \left(\partial{g}(f(v))\partial f(v)\right)$,
which shows~\eqref{cr-2},
when $\ell=1$. The proof for $\ell>1$ follows readily
from the result when $\ell=1$, together with the observation
that the convex hull of a rectangular set is rectangular.
\end{proof}
We  get as an easy corollary 
one inclusion in Theorem~\ref{ncf}:
\begin{cor}\label{ncf2}
If $f$ is as in Theorem~\rm{\ref{ncf}}, and if $k\geq 1$,
then, $\partial{f}^k(v)=\vex\big((\partial{f}(v))^k\big)$,
for all eigenvectors $v$ of $f$.
Moreover, $\gc(f^k)\supset \gc(f)^k$.
\end{cor}
\begin{proof}
By~\eqref{cr-2}, $\partial{f}^k(v)=\vex (\partial{f(v)}\vex(\cdots
\vex \partial{f(v)}\cdots))=\vex (\partial{f(v)}\cdots \partial{f}(v))
=\vex \big((\partial{f}(v))^k\big)$.
Assertion~\ref{b2} of Corollary~\ref{crit}
shows that there is a matrix $P\in \partial{f}(v)$
such that $\gc(f)=\gf(P)$. Since $P^k\in (\partial{f}(v))^k\subset
\partial{f^k}(v)$,
$\gf(P^k)\subset \gc(f^k)$, and, using Proposition~\ref{cor-c}, we get
$\gc(f)^k=\gf(P)^k=\gf(P^k)\subset \gc(f^k)$.
\end{proof}
To show the other inclusion, we
shall need the following nonlinear version of Proposition~\ref{propf},
which relates final classes with invariant measures.
\begin{prop}\label{mesf}
Let $f$ denote a convex monotone homogeneous map $\R^n\to\R^n$
that has a fixed point, $v$, and let $F$ be a subset  of $\{1,\ldots,n\}$.
The following propositions are equivalent:
\begin{enumerate}
\item  $F$ is a union of elements of $\cc(f)$,
\item  there exists a stochastic  row  vector $m$ with support $F$,
such that $mf(x)\geq mx,\;\forall x\in\R^n$,
\item there exists a stochastic  row  vector $m$ with support $F$,
such that $mf'_v(x)\geq mx,\;\forall x\in\R^n$.
\end{enumerate}
\end{prop}
\begin{proof}
Without loss of generality, we assume that $v=0$ (otherwise,
replace $f$ by $-v+f(\cdot +v)$).

\noindent $(1 \implies 2)$
Let $F$ be a union of elements of $\cc(f)=\cf(\partial f(0))$.
Using the first assertion of
Proposition~\ref{tech}, we can write $F$ as the union of
disjoint sets $F_1,\ldots , F_k\in\cc(f)$.
Moreover, there is a matrix $P\in\partial f(0)$
such that $F_1,\ldots , F_k$ are final classes of $P$.
Indeed, pick $P^{(1)},
\ldots,P^{(k)}\in \partial f(0)$
such that for all $1\leq l\leq k$,
$F_l$ is a final class of $P^{(l)}$.
We denote by $P$ the matrix whose $i$-th
row is equal to the $i$-th row of $P^{(l)}$, when  $i\in F_l$, 
and to any element of $\partial f_i(0)$, when
$i\not\in F$.
By rectangularity of $\partial f(0)$, $P\in \partial f(0)$, and
$F_1,\ldots, F_k$ are final classes of $P$.

Thus, $F$ is a union of final classes of $P$.
It follows, by Proposition~\ref{propf}, that $P$ has an invariant measure $m$
with support $F$. Since $P\in\partial f (0)$, we get $f(x)\geq Px$
for all $x\in\R^n$, so $mf(x)\geq mPx=mx$.

\noindent  $(2\implies 3)$ 
Assume that there exists  
a stochastic  row  vector $m$ with support $F$,
such that $mf(x)\geq mx$ for all $x\in\R^n$.
Then,  $m \varepsilon^{-1} f(\varepsilon x)\geq mx$, for all $\varepsilon>0$.
Taking the limit when $\varepsilon$ goes to zero,
one gets $mf'_0(x)\geq mx$ for all $x\in\R^n$.

\noindent $(3 \implies 1)$
Assume now that there exists 
a stochastic  row  vector $m$ with support $F$,
such that $mf'_0(x)\geq mx$ for all $x\in\R^n$.
Applying Property \eqref{propfp} to the coordinates of $f'_0$ 
and using the rectangularity of $\partial f(0)$,
we get
\[mf'_0(x)=m\sup_{P\in \partial f(0)} Px = \sup_{P\in \partial f(0)} mPx
\enspace.\]
Hence, 
setting $\sM=\{ mP,\;  P\in \partial f(0)\}$,
we see that the assumption on $m$ is equivalent to
 $\sup_{\nu\in \sM} \nu x\geq mx$ for all $x\in \R^n$.
In order to prove point 1, it is sufficient to
show that there exists $P\in\partial f(0)$ such that
$mP=m$, that is $m\in\sM$.
 Indeed, by Proposition~\ref{propf}, this would imply that $F$ is a union
 of final classes of $P$, thus a union of elements of $\cc(f)$.
Assume, by contradiction, that $m\not\in \sM$.
Since  $\partial f(0)$ is a convex compact subset of $\R^{n\times n}$,
$\sM$ is a convex compact subset of the space $\R^{1 \times n}$
of row vectors. Hence, the Hahn-Banach theorem implies that
there exists a separating hyperplane, that is there exists a column vector
$x\in \R^n$ and  a real $\alpha$, such that
$mx>\alpha$ and $\nu x<\alpha$ for all $\nu\in\sM$.
Then,  $\sup_{\nu\in \sM} \nu x\leq \alpha<mx$, a contradiction.
We have proved that $m\in \sM$, which shows point 1.
\end{proof}
We shall also need the following notion of \new{additive recession function}.
Let $g:\R^n\to\R$ be an additively \new{subhomogeneous map},
that is a map such that $g(\lambda +x)\leq \lambda +g(x)$ for
all $\lambda\geq 0$ and $x\in\R^n$.
The sequence  $-\rho+g(\rho+y)$ is nonincreasing with respect
to $\rho>0$, allowing us to define
the additive recession function $\tilde{g}:\R^n\to \Rmax$ of $g$ by
\begin{equation}\label{defghat2}
 \tilde{g}(y)=\lim_{\rho\to+\infty} -\rho+g(\rho +y) \enspace. 
\end{equation}
If $g$ is monotone and convex, so does $\tilde{g}$.
Moreover, the domain of $g^*$ is included in the set 
$\probam{n}=\set{p\in \R^n}{\sum_{1\leq i\leq n}p_i\leq 1,\;
p_1,\ldots,p_n\geq 0}$ of \new{substochastic} vectors, and $\tilde{g}$
is given by:
\begin{equation}\label{defghat}
 \tilde{g}(y)= \sup_{q\in\dom g^*\cap \proba{n}} q \cdot y-g^*(q)\enspace,
 \end{equation}
where an empty supremum is equal to $-\infty$.
By \eqref{defghat2} and the nonexpansiveness of $g$, $\tilde{g}$
is either finite everywhere or identically equal to $-\infty$.

We define similarly the additive recession function of an
additively subhomogeneous map $g:\R^n\to\R^n$. If $g=(g_1,\ldots,g_n)$, 
we have $\tilde{g}=(\widetilde{g_1},\ldots,\widetilde{g_n})$.
If $g$ is monotone and convex, and if $\tilde{g}$ has only finite entries,
we have,
by \eqref{defghat2},
$g(\rho+y)=\rho+\tilde{g}(y)+o(1)$, when $\rho$ goes to $+\infty$. 
So using the nonexpansiveness of $g$, we get 
$g\circ h (\rho+y)=g(\rho+\tilde{h}(y))+o(1)=\rho+\tilde{g}\circ
 \tilde{h}(y)+o(1)$,
for any other map $h$ with the same properties as $g$, 
which leads to the following chain rule:
\begin{lemma}\label{chainhat}
Let $g$ and $h:\R^n\to\R^n$ be two monotone, convex and subhomogeneous maps,
such that $\tilde{g}$ and $\tilde{h}$ take finite values.
Then, $g\circ h:\R^n\to\R^n$ is monotone, convex, subhomogeneous,
and $\widetilde{g\circ h}=\tilde{g}\circ \tilde{h}$.
\end{lemma}

For any $n\times n$ stochastic matrix $P$, we say that the set
 $F\subset \{1,\ldots, n\}$
is \new{invariant} by $P$ if $P_{ij}=0$ for all $i\in F$ and $j\not \in F$
(that is, if it is invariant by the
dynamics of the Markov chain of transition matrix $P$).

\begin{lemma}\label{univ}
Let $f$ denote a convex monotone homogeneous map $\R^n\to\R^n$
that has an eigenvector $v$. Then, there exists at least one
critical class of $f$ that is invariant by all
the matrices $P\in\partial f(v)$.
We say that such a critical class is \new{invariant} by $f$.
\end{lemma}
\begin{proof}
Let $\sG=\cup_{P\in\partial f(v)}\sG(P)$ and 
$F$ denotes an
arbitrary final class of $\sG$.
The set $F$ is clearly invariant by all the matrices $P\in \partial f(v)$.
Since $\sG$ is finite, there exists a finite set 
$\sP\subset \partial{f}(v)$
such that $\sG=\cup_{P\in\sP}\sG(P)$.
Consider $\widehat{P}=|\sP|^{-1}(\sum_{P\in \sP} P)$.
Using the convexity of $\partial{f}(v)$, 
we have $\widehat{P}\in \partial{f}(v)$.
Since the graph of $\widehat{P}$ is equal to $\sG$,
$F$ is a final class of $\widehat{P}$, hence $F\in \cc(f)$.
This is a critical class of $f$ since for all $P\in \partial{f}(v)$,
$P_{ij}=0$ for all $i\in F$ and $j\not\in F$.
\end{proof}
\begin{example}
\label{ex-f4}
For the map $f:\R^3\to \R^3$,
\[
f(x) = \left(\begin{array}{c}
x_2\\
x_2 \vee x_3\\
x_3 
\end{array}
\right)
\enspace,
\] 
both $\{2\}$ and $\{3\}$ are critical
classes, but only $\{3\}$ is invariant.
\end{example}
The last ingredient in the proof of Theorem~\ref{ncf}
is a simple restriction operation. If $N$ is a subset of $\{1,\ldots,n\}$,
we denote by $\sfi_N$ the trivial injection
$\R^N\to \R^n$, which is defined by $(\sfi_N(x))_i=x_i$
if $i\in N$, and $(\sfi_N(x))_i=0$, otherwise.
We define the restriction $f_{NN}:\R^N\to \R^N$ of $f$ to $N$ by
$f_{NN}= \sfr_N\circ f\circ \sfi_N$.
If $f$ is convex monotone and bihomogeneous, then 
$f_{NN}$ is clearly monotone, convex, multiplicatively
homogeneous, and additively subhomogeneous.
Indeed, $f(x)=f'_0(x)=\sup_{P\in \partial f(0)} Px$,
so $f_{NN}(y)=\sup_{Q\in \partial{f_{NN}}(0)}Qy$,
with $\partial{f_{NN}}(0)=\set{P_{NN}}{P\in \partial f(0)}$.
Recall that the \new{critical graph} of
an additively \new{subhomogeneous} map has
been defined in~\S\ref{subsec-subhom}.
We have the following lemma, which follows readily
from the fact that
for any $P\in\partial f(0)$, $P_{FN}=0$, so that any final class of $P$
is included either in $N$ or in $F$.
\begin{lemma}\label{lem-restriction}
Let $F$ denote an invariant critical class of a convex monotone 
bihomogeneous map $f:\R^n\to \R^n$,
and let $N=\{1,\ldots,n\}\setminus F$. Then, 
$f\circ \sfi_N = \sfi_N \circ f_{NN}$, $f_{FF}\circ \sfr_F = \sfr_F \circ f$
and the map $f_{FF}$ is additively homogeneous.
Moreover, $\gc(f)=\gc(f_{NN})\cup \gc(f_{FF})$, 
and, when all the entries of $\widetilde{f_{NN}}$ are finite,
\begin{equation}
\gc(f)=\gc(\widetilde{f_{NN}})\cup \gc(f_{FF})
\label{e-disjointu}
\end{equation}
(the unions are disjoint).\qed
\end{lemma}
In general, $f_{NN}$ is only additively subhomogeneous,
 since the $N\times N$ submatrices $P_{NN}$ are
only substochastic.
\begin{example}
For the map $f$ of Example~\ref{ex-f4},
we have $F=\{3\}$, $N=\{1,2\}$, and
\[
f_{NN} 
\left(\begin{array}{c}
x_1\\
x_2 
\end{array}
\right)
= \left(\begin{array}{c}
x_2\\
x_2 \vee 0
\end{array}
\right)
\enspace,
\quad 
\widetilde{f_{NN}}
\left(\begin{array}{c}
x_1\\
x_2 
\end{array}
\right)
= \left(\begin{array}{c}
x_2\\
x_2 
\end{array}
\right)
\enspace,
\quad 
f_{FF}(x_3)=x_3 \enspace ,
\] 
and $\gc(f)=\gc(\widetilde{f}_{NN})\cup \gc(f_{FF})$
is the disjoint union of the two graphs
consisting of a single loop at node $2$, and $3$,
respectively.
\end{example}
\begin{proof}[Proof of Theorem~\ref{ncf}]
Without loss of generality, we assume that $0$ is
a fixed point of $f$.
By definition of $\gc(f)$, $\gc(f)=\gc(f'_0)$.
Moreover,  Lemma~\ref{chain} shows that $(f^k)'_0 = (f'_0)^k$ so that
 $\gc(f^k)=\gc((f^k)'_0)= \gc((f'_0)^k)$.
Thus, it is enough to prove Theorem~\ref{ncf} for maps of the form
$f'_0$, that is for convex monotone homogeneous maps that are
also multiplicatively homogeneous.

We next reduce to the case where $\nc(f^k)=\{1,\ldots , n\}$.
Indeed, let $N= \{1,\ldots , n\}\setminus \nc(f^k)$ and consider
$g:\R^n\to\R^n$ defined by $g_i(x)=f_i(x)\vee x_i$ if $i\in N$ and
$g_i(x)=f_i(x)$ otherwise.
Using the inclusion already proved in Corollary~\ref{ncf2},
we get $\nc(f)\subset \nc(f^k) = \{1,\ldots , n\}\setminus N$, hence
$\gc(g)$ is the disjoint union of $\gc(f)$ and of the set $\sL$ of
loops $i\to i$ with $i\in N$. This implies in particular that
$\gc(g)^k$ is the disjoint union of  $\gc(f)^k$ and $\sL$,
and that  $\nc(g)=\nc(f)\cup N$.
Since $g(0)=0$ and
$g^k\geq f^k$, we get $\nc(g^k)\supset \nc(f^k)$ and $\gc(g^k)\supset \gc(f^k)$.
Again by Corollary~\ref{ncf2},
we get $\nc(g^k)\supset \nc(g)\supset N$, hence,
$\nc(g^k)=\{1,\ldots , n\}$.
Now if $\gc(g^k)\subset\gc(g)^k$, then $\gc(f^k)\subset \gc(g^k)\subset
\gc(g)^k=\gc(f)^k\cup \sL$, and since $\nc(f^k)= \{1,\ldots , n\}\setminus N$,
we have $\gc(f^k)\subset \gc(f)^k$, which, together
with the inclusion proved in 
Corollary~\ref{ncf2}, yields
$\gc(f^k)= \gc(f)^k$.

We are now reduced to show the following proposition,
$(\mathsf{P}_n)$:
{\em for all convex monotone bihomogeneous
maps $f:\R^n\to\R^n$, and for all $k\geq 1$, such that
$\nc(f^k)=\{1,\ldots, n\}$, we have
$\gc(f^k)\subset\gc(f)^k$}.
We shall prove $(\mathsf{P}_n)$ by induction
on the dimension $n$.
Proposition $(\mathsf{P}_1)$ trivially holds, since in dimension $1$,
the only possible map $f$ is the identity map.
Assume now that $(\mathsf{P}_\ell)$ holds for all dimensions 
$\ell<n$, and consider a map $f:\R^n\to\R^n$ as in $(\mathsf{P}_n)$.

Let $F$ be an invariant critical class of $f$ (as in Lemma~\ref{univ}),
let $N=\{1,\ldots, n\}\setminus F$, and
put $g=f_{NN}$ and $h=f_{FF}$.
Assume first that all the entries of $\widetilde{g}$ are finite.
Then, by~\eqref{e-disjointu}, $\gc(f)=\gc(\widetilde{g})\cup \gc(h)$.
Since $F$ is also invariant for $f^k$, we also have 
$\gc(f^k)=\gc(\widetilde{(f^k)_{NN}})\cup \gc((f^k)_{FF})$.
Using the first part of Lemma~\ref{lem-restriction}, we obtain 
$(f^k)_{NN}= \sfr_N\circ f^k\circ i_N=g^k$ and $(f^k)_{FF}=h^k$.
Moreover, if the entries of $\widetilde{g}$ are finite,
Lemma~\ref{chainhat} yields
$\widetilde{g^k}=\widetilde{g}^k$,
so $\gc(f^k)=\gc(\widetilde{g}^k)\cup \gc(h^k)$.
Thus, in order to prove $\gc(f^k)\subset \gc(f)^k$, it is
sufficient to prove the same inclusion for $h$ and $\widetilde{g}$
and to prove that $\widetilde{g}$ has finite entries.

Let us first prove  $\gc(h^k)\subset \gc(h)^k$.
Since $F$ is an
invariant critical class of $f$,
$\gc(h)$ is strongly connected and contains all the nodes of $F$.
By Assertion~\ref{b2} of Corollary~\ref{crit}, there exists 
$P\in \partial h(0)$
such that $\gc(h)=\gf(P)$, hence $P$ is irreducible and $\gc(h)=\sG(P)$.
Moreover, for all $Q\in\partial h(0)$, the matrix
$\frac{Q+P}{2}$
belongs to $\partial h(0)$
and is irreducible, so $\gf(Q)\subset\sG(Q)\subset \sG(\frac{Q+P}{2})
=\gf(\frac{Q+P}{2})\subset \gc(h)$. 
Defining $\sG(\sP)=\cup_{P\in\sP} \sG(P)$
for all sets $\sP$ of stochastic matrices,
we get $\gc(h)=\sG(\partial h(0))$.
Hence, $\gc(h^k)\subset \sG(\partial (h^k) (0))=\sG(\vex ((\partial h(0))^k))
=\sG((\partial h(0))^k)\subset \sG(\partial h(0))^k=\gc(h)^k$.

Let us now prove that all the entries of $\widetilde{g}$ are finite.
This is indeed equivalent to the property $g(1)=1$, where
$1$ denotes the unit vector (which means that
$ \partial g(0)\cap \proba{NN}\not =\varnothing$).
We know that $g(1)\leq 1$, so by monotonicity of $g$,
$g^k(1)\leq \cdots\leq g(1)\leq 1$.
On the other hand, by Proposition~\ref{mesf}, there exists 
a stochastic (row) vector $m$ with positive entries (i.e., the support 
of $m$ is $\{1,\ldots n\}$), such that $mf^k(x)\geq mx$ for all $x\in\R^n$.
By Lemma~\ref{lem-restriction},  $f\circ \sfi_N=\sfi_N \circ g$,
hence $f^k\circ \sfi_N=\sfi_N \circ g^k$.
Applying  $mf^k(x)\geq mx$ to $x= \sfi_N(y)$, we
obtain  $m_Ng^k(y)\geq m_Ny$, where $m_N$ is the restriction of $m$ to the
set $N$ ($m_N$ has positive entries since all the entries
of $m$ are positive).  
So $m_N 1 \geq m_N g(1)\geq \cdots \geq m_N g^k(1)\geq m_N 1$.
This yields $m_N 1 = m_N g(1)$. Using the positivity of $m_N$
and the inequality $g(1)\leq 1$, we get $1=g(1)$.

We finally show that $\gc(\widetilde{g}^k)\subset \gc(\widetilde{g})^k$,
by using the induction assumption.
Indeed, since $F\not=\varnothing$, the cardinality of $N$, which
is the dimension of the space
on which $g$ operates, is strictly less than $n$.
Moreover, since $m_Ng^k(y)\geq m_Ny$ for all $y\in \R^N$,
we deduce that $m_N\widetilde{g^k}(y)\geq m_Ny$.
Since  $\widetilde{g}$ has finite entries, $\widetilde{g^k}=\widetilde{g}^k$,
so that $m_N\widetilde{g}^k(y)\geq m_Ny$ for all $y\in \R^N$, and
Proposition~\ref{mesf} shows that $\nc(\widetilde{g}^k)=N$.
The induction assumption can then be applied to $\widetilde{g}$, 
and yields $\gc(\widetilde{g}^k)\subset \gc(\widetilde{g})^k$.
This concludes the proof of $(\mathsf{P}_n)$.
\end{proof}

\section{Cyclicity Theorem for Convex Monotone Homogeneous Maps}
\label{cyclicity}
In this section, we use our knowledge of the eigenspace
of $f$ to study the asymptotic behavior of $f^k$ when
$k$ tends to $+\infty$.  In particular, we are
interested in the \new{periodic orbits}
of $f$, which are of the form $\{f^k(x)\}_{k\in \N}$,
with $f^c(x)=x$ for some $c\geq 1$. 
The set of such $c$ is exactly the set of multiples of a positive integer,
which is the \new{length} of the orbit $\{f^k(x)\}_{k\in \N}$.

Let us first recall some more or less classical
facts on periodic orbits of stochastic matrices.
The \new{cyclicity} $\sfc(\sG)$ of a strongly connected
graph $\sG$ is the $\gcd$ of the lengths of the circuits of $\sG$.
The \new{cyclicity} of a graph $\sG$ with strongly connected
component $\sG_1,\ldots, \sG_s$ is defined by $\sfc(\sG)=\lcm(\sfc(\sG_1),
\ldots,\sfc(\sG_s))$. For a stochastic matrix $P$,
we define the \new{cyclicity} of $P$ as: $\sfc(P)=\sfc(\gf(P))$.
The name ``cyclicity'' is justified by the following result.
\begin{prop}\label{prop-class}
If $P$ is a stochastic matrix, then, the length of
any periodic orbit of $x\mapsto Px$ divides $\sfc(P)$.
\end{prop}
\begin{proof}
If $P$ is a stochastic matrix with cyclicity $c=\sfc(P)$,
then the Perron-Frobenius theorem shows that $P^{\ell c}\to \Pi$
when $\ell\to\infty$, where $\Pi$ is the
spectral projector of $P^c$ for the eigenvalue
$1$. 
If $P^k x=x$, then, $x=\lim_{\ell\to\infty} P^{\lcm(k,c)\ell}x=\Pi x$.
Since the spectral projector satisfies $P^c\Pi =\Pi$,
we get $x=P^c x$, which shows that
the length of any periodic
orbit of $P$ divides $c$.
\end{proof}
More generally, we define the \new{cyclicity}
$\sfc(f)$ of a convex monotone homogeneous map $f$
as the cyclicity of its critical graph $\sfc(\gc(f))$.

\begin{prop}
\label{attained}
Let $f$ denote a convex monotone homogeneous map
that has an eigenvector $v$. Then, 
\begin{equation}
\sfc(f)= \gcd\set{\sfc(P)}{P\in \partial f(v),\; \nf(P)=\nc(f)}  \enspace .
\label{eq-def-cyc}
\end{equation}
Moreover, there exists $P\in \partial{f}(v)$ such that
$\nf(P)=\nc(f)$ and $\sfc(P)=\sfc(f)$.
\end{prop}
\begin{proof}
Let $P\in \partial{f}(v)$ such that $\nf(P)=\nc(f)$.
Since $\gf(P)\subset \gc(f)$, and since both
graphs have the same set of nodes,
each critical class $C$ of $f$ can be partitioned into
final classes of $P$, $C=F_1\cup\ldots\cup F_r$,
and $\gf(P)|_C\subset \gc(f)|_C$, where $\sG|_C$ denotes the restriction
of the graph $\sG$ to $C$.
Hence, the cyclicity of $\gc(f)|_C$, which is equal
to the $\gcd$ of the lengths of circuits of $\gc(f)|_C$ divides
 $\sfc(\gf(P)|_C)=\lcm(\sfc(P_{F_1F_1}),\ldots,\sfc(P_{F_r,F_r}))$.
Since this divisibility holds for all critical
classes $C$ of $f$, $\sfc(f)$ divides $\sfc(P)$.
Conversely, Assertion~\ref{b2} of Corollary~\ref{crit},
there is a matrix $P\in\partial{f}(v)$
such that $\gf(P)=\gc(f)$, which shows that the right
hand side of~\eqref{eq-def-cyc}
divides $\sfc(f)$.
\end{proof}
\begin{prop}\label{orbitclass}
Let $f$ be a convex monotone homogeneous map that has a 
fixed point  and let $\sG_1,\ldots,\sG_s$ be
the strongly connected components of $\gc(f)$.
Then, $f^{\sfc(f)}$ has $\sfc(\sG_1)+\cdots+\sfc(\sG_s)$ 
critical classes and has cyclicity~$1$.
\end{prop}
\begin{proof}
From Theorem~\ref{ncf},
we get $\gc(f^{\sfc(f)})=\gc(f)^{\sfc(f)}=\cup_{i=1}^s\sG_i^{\sfc(f)}$
(disjoint union).
For any strongly connected graph $\sG$
with cyclicity $c$ and any integer $k$, $\sG^{kc}$ has
$c$ strongly connected components with cyclicity~$1$ 
(this is a well known result which follows from the decomposition
of $\sG$ into cyclic classes, see e.g. \cite[Chapter 6, \S3]{kemeny} 
and \cite[Chapter 2, \S2]{berman}).
Since $\sfc(\sG_i)$ divides $\sfc(f)$,
we obtain that $\gc(f^{\sfc(f)})$ has  $\sfc(\sG_1)+\cdots+\sfc(\sG_s)$  
strongly connected components with cyclicity~$1$.
\end{proof}
We shall also need the following classical fact.
\begin{prop}\label{exists}
A monotone homogeneous map $\R^n\to\R^n$
which has a periodic orbit has a fixed point.
\end{prop}
\begin{proof}
If $f$ has a periodic orbit,
there is an $x\in \R^n$ and $k\geq 1$
such that $f^k(x)=x$.
Consider $z=x\wedge f(x)\wedge \ldots \wedge f^{k-1}(x)$,
By monotonicity of $f$, $f(z)\leq f(x)\wedge f^2(x)\wedge \ldots \wedge f^k(x)
=z$. Thus, the sequence $\{f^\ell(z)\}_{\ell\in \N}$ is nonincreasing,
and since by nonexpansiveness of $f$,
$|f^\ell(z)-f^\ell(x)|\leq |z-x|$ is bounded, 
$f^\ell(z)$ which is bounded and nonincreasing converges
to a limit  which is an eigenvector of $f$.
\end{proof}
\begin{theorem}\label{theo-cycli}
The length of any periodic orbit
of a convex monotone homogeneous map $f$
divides $\sfc(f)$.
\end{theorem}
\begin{proof}
If $f$ has a periodic orbit $\{f^k(x)\}_{k\in \N}$,
then, by Proposition~\ref{exists},
$f$ has a fixed point $v$.
It suffices to prove Theorem~\ref{theo-cycli} when $\sfc(f)=1$.
Indeed, if  $\sfc(f)\not =1$, consider $g=f^{\sfc(f)}$, which is
such that  $\sfc(g)=1$ (by Proposition~\ref{orbitclass}) and
 $\{g^k(x)\}_{k\in \N}$ is periodic.
If we know that Theorem~\ref{theo-cycli} holds for $g$, we get
$x=g(x)=f^{\sfc(f)}(x)$, hence
the length of the periodic orbit $\{f^k(x)\}_{k\in \N}$ divides $\sfc(f)$,
which yields Theorem~\ref{theo-cycli} in the general case.

Let us now prove the theorem when $\sfc(f)=1$.
Using Proposition~\ref{attained}, we can find a matrix $P\in \partial{f}(v)$
such that $\sfc(P)=\sfc(f)=1$ and $\nf(P)=\nc(f)$. Using $f(v)=v$
and $f(y)-f(v)\geq P(y-v)$, we get after an immediate
induction, 
\begin{equation}
\label{induc}
f^k(x)\geq P^k(x-v)+ v \enspace \forall k\in \N \enspace.
\end{equation}
In particular, 
if  $c$ is the length of the orbit $\{f^k(x)\}_{k\in \N}$,
$x=f^c(x)\geq P^c(x-v)+v$.

To make the proof more intuitive, we
shall first show the theorem when $\nc(f)=\{1,\ldots,n\}$.
By Proposition~\ref{cor-c},
$P^c$ again has $\{1,\ldots,n\}$ as union of final classes,
and applying Lemma~\ref{disc} to $x-v\geq P^c(x-v)$,
we get $x-v=P^c(x-v)$. Since $P$ has cyclicity~$1$, this implies that
$P(x-v)=x-v$, and using~\eqref{induc},
we get $f(x)\geq x$. By monotonicity of $f$, we get 
$x=f^c(x)\geq \cdots \geq f(x)\geq x$.
Thus, $f(x)=x$, and $c=1=\sfc(f)$.

Consider now the general case where $F=\nc(f)\not=\{1,\ldots,n\}$.
As before, we get 
$x-v\geq P^c (x-v)$. By Lemma~\ref{disc} and 
Proposition~\ref{cor-c},
this yields $x-v= P^c (x-v)$ on $\nf(P^c)=\nf(P)=F$. Since $P_{FN}=0$
for $N=\{1,\ldots, n\}\setminus F$, this equation can be rewritten as
$(x-v)_F= (P_{FF})^c(x-v)_F$. Since $c(P_{FF})=c(P)=c(f)=1$, we deduce
$(x-v)_F= P_{FF}(x-v)_F$.
Then, for all $k\geq 1$, $(P_{FF})^{k}(x-v)_F=(x-v)_F$
or equivalently  $(P^{k}(x-v))_F=(x-v)_F$.
Putting this in \eqref{induc}, we get $f^{k}(x)\geq x$ on $F$.
Consider now $z=x\wedge f(x)\wedge\cdots\wedge f^{c-1}(x)$,
which coincides with $x$ on $F$.
By monotonicity of $f$, $f(z)\leq z$.
Hence, by Lemma~\ref{lem-projection}, $w=f^\omega(z)$ is a fixed
point of $f$ and it coincides with $z$ on $\nc(f)=F$.
Since $f(w)=w$, we get $f^c(w)=w$.
So $w$ and $x$ are two fixed points of $f^c$, which coincide on
$F=\nc(f)=\nc(f^c)$. This implies that $w=x$ by Theorem~\ref{thmain}.
Hence, $f(x)=x$, and by the minimality of $c$ (as the
length of the orbit of $x$), we get $c=1$.
\end{proof}
\begin{cor}
The lengths of periodic orbits of convex
monotone homogeneous maps $\R^n\to\R^n$ 
are exactly the orders of the elements of the symmetric
group on $n$ letters.\qed
\end{cor}
Combining Theorem~\ref{theo-cycli} with the theorem
of Nussbaum~\cite{nuss90} and Sine~\cite{sine}, we get:
\begin{cor}\label{cor-asymp}
If $f$ is a convex monotone homogeneous
map $\R^n\to \R^n$ with eigenvalue $\lambda$,
then, for all $x\in \R^n$, $f^{k\sfc(f)}(x)-k\sfc(f)\lambda$
converges when $k\to \infty$.\qed
\end{cor} 
In other words, introducing the normalized
map $g=f-\lambda$, which is such that $g^k=f^k-k\lambda$,
we see that any orbit $\{g^k(x)\}_{k\in \N}$ of
$g$ converges to a periodic orbit of length at most $\sfc(g)=\sfc(f)$.
We also get, as an immediate consequence of Theorem~\ref{theo-cycli}:
\begin{cor}\label{c-cycdim}
For any convex monotone homogeneous map $f$ that has a 
fixed point, the union of the periodic orbits of $f$
is exactly the set $\sE(f^{\sfc(f)})$ of fixed points of $f^{\sfc(f)}$.
\qed
\end{cor}
Combining Corollary~\ref{c-cycdim}, Proposition~\ref{orbitclass},
Theorem~\ref{ncf} 
and Theorem~\ref{thmain}, we obtain:
\begin{cor}\label{c-cycdim2}
Let $f$, $\sG_1,\ldots,\sG_s$ be as in Proposition~\ref{orbitclass},
and let $C=\nc(f)$.
Then, $\sfr_C$ sends bijectively the union of 
periodic orbits, $\sE(f^{\sfc(f)})$, to a convex inf-subsemilattice of $\R^C$, 
whose dimension 
is at most equal to $\sfm(f^{\sfc(f)})=\sfc(\sG_1)+\cdots+\sfc(\sG_s)$.\qed
\end{cor}

\section{Piecewise Affine Convex Monotone Homogeneous Maps}
\label{piecewise}
\subsection{Dimension of the Eigenspace}
As discussed in \S\ref{sec-rel}, the dimension of the eigenspace
of a piecewise affine 
convex monotone homogeneous map was characterized by Romanovsky~\cite{roma73}
and by Schweitzer and Federgruen~\cite{sch78}: this 
shows that the bound on the dimension
given in Corollary~\ref{cor-summarize}, is attained when
$f$ is piecewise affine. 
In this subsection, we give an independent proof
of this fact, which shows some
qualitative properties of $\sE(f)$ (connection
between $\sE(f)$ and $\sE(f'_v)$, for any
eigenvector $v$, role of invariant critical
classes).

Let us first recall some basic
definitions and facts. A \new{polyhedron} is an
intersection of finitely many half-spaces.
A map $f:\R^n\to \R$ is \new{piecewise affine} if
$\R^n$ can be covered by finitely many polyhedra such that
the restriction of $f$ to each polyhedron is affine.
The following result is an immediate consequence
of classical results on convex maps with polyhedral
epigraphs~\cite{ROCK}.
\begin{prop}\label{prop-zero}
A convex map $f:\R^n\to \R$ is piecewise affine if, and
only if, there exists a finite set $\sP\subset \dom f^*$
such that 
\begin{equation}
f(x)= \max_{p\in \sP}(p\cdot x - f^*(p)) \enspace.
\label{sp}
\end{equation}
\end{prop}
It will be useful to consider bihomogeneous maps.
\begin{cor}\label{cor-canon}
A piecewise affine convex monotone bihomogeneous
map $f:\R^n\to \R$ can be written as
\begin{equation}
f(x)=\max_{p\in \sP}p\cdot x \enspace,
\label{sp2}
\end{equation}
where $\sP$ is a finite set of stochastic vectors.
\end{cor}
\begin{proof}
By Proposition~\ref{prop-0},
$\dom f^*\subset \proba n$,
therefore, the finite set $\sP$ in~\eqref{sp}
is composed of stochastic vectors,
and since $f^*$ takes only the values
$0$ and $+\infty$ when $f$ is multiplicatively
homogeneous, we get~\eqref{sp2}.
\end{proof}
The following simple fact
shows that eigenspaces of convex piecewise affine
maps can be effectively computed.
\begin{cor}\label{cor-compute}
The eigenspace of a convex piecewise affine monotone homogeneous
map $f:\R^n\to \R^n$ is a finite union
of polyhedra.
\end{cor}
\begin{proof}
We assume, without loss of generality, that $f$ has eigenvalue $0$.
By Propositions~\ref{prop-0} and~\ref{prop-zero},
we can write each coordinate of $f$
as $f_i(x)=\max_{p\in \sP_i} (p\cdot x -f_i^*(p))$,
where $\sP_i$ is a finite subset of $\proba n$.
To any $i\in \{1,\ldots,n\}$ and $p\in \sP_i$,
we associate the (possibly empty) polyhedron
$K_{i,p}=\set{x\in \R^n}{x_i=p\cdot x-f_i^*(p)
\geq q\cdot x -f_i^*(q),\forall q\in \sP_i}$.
We denote by $\Phi$ the set
of maps $\varphi:\{1,\ldots,n\}
\to \sP_1\cup\ldots\cup \sP_n$ such that
$\varphi(i)\in \sP_i$, for all $i\in \{1,\ldots,n\}$.
To $\varphi\in \Phi$, we associate
the polyhedron $K_{\varphi}=\cap_{1\leq i\leq n}K_{i,\varphi(i)}$.
If $x=f(x)$, then, for each $i$, by finiteness of $\sP_i$,
we have $x_i=\max_{p\in \sP_i} p\cdot x-f_i^*(p)=
\varphi(i) \cdot x-f_i^*(\varphi(i))$,
for some $\varphi(i)\in\sP_i$. This shows
that $\sE(f)\subset \cup_{\varphi\in \Phi}K_\varphi$.
Since the other inclusion holds trivially,
$\sE(f)= \cup_{\varphi\in \Phi}K_\varphi$
is a finite union of polyhedra.
\end{proof}
The proof of Corollary~\ref{cor-compute} yields
an algorithm (with exponential execution time)
to compute $\sE(f)$, but it does not
tell much about the geometry of $\sE(f)$. We shall next
prove more qualitative properties.

The following observation
shows that the directional derivative $f'_v$  defined in 
Section~\ref{secttang} is a ``tangent'' map of $f$.
\begin{lemma}\label{lemmafv}
If $f:\R^n\to \R$ is a piecewise affine convex map,
then, for all $v\in \R^n$, there is a neighborhood $V$
of $v$ such that 
\begin{equation}
\label{neighbor}
f(x)= f(v)+f'_v(x-v)\qquad \forall x\in V
\enspace .
\end{equation}
\end{lemma}
\begin{proof}
The definition of $\partial f(v)$
together with \eqref{propfp} show that the inequality
$\geq$ holds in~\eqref{neighbor}. To show the
other inequality, we write $f$ as~\eqref{sp},
and set $\sQ=\sP\cap \partial f(v)$,
$\sR=\sP\setminus \sQ$.
We can write $f=g\vee h$ where 
\begin{equation}
g(x)=
\max_{p\in \sP\cap \partial f(v)}
(p\cdot x-f^*(p))
= f(v)+\max_{p\in \sQ}
p\cdot (x-v)
\label{sp4}
\end{equation}
(since $f(v)=p\cdot v- f^*(p)$,
for all $p\in \partial f(v)$),
and 
$h(x)=\max_{p\in \sR} (p\cdot x-f^*(p))$.
The functions $g$ and $h$ are continuous
and satisfy $f(v)=g(v)>h(v)=\max_{p\in \sP\setminus\partial f(v)}
(p\cdot v-f^*(v))$, where the strict inequality follows
from the characterization~\eqref{e-diesis}
of $\partial f(v)$ and
from the finiteness of $\sP$. Therefore,
there is a neighborhood $V$ of $v$ such that $f(x)=g(x)\vee h(x)=g(x)$
for all $x\in V$. Combining this with~\eqref{sp4} and \eqref{propfp},
we get $f(x)=g(x)\leq f(v)+f'_v(x-v)$ for all $x\in V$.
\end{proof}
Our interest in $f'_v$ stems from the following reduction.
\begin{lemma}\label{reduc1}
Let $f$ be 
a piecewise affine convex monotone homogeneous map
$\R^n\to\R^n$ with eigenvector $v$.
Then, the dimensions of $\sE(f)$
and $\sE(f'_v)$ are the same.
\end{lemma}
\begin{proof}
We first note that by construction of $f'_v$,
the critical classes for $f$ and $f'_v$ are
the same. Let $C=\nc(f)=\nc(f'_v)$.
By Corollary~\ref{cor-summarize},
there is a set $U\subset \sE(f'_v)$
such that $\sfr_C(U)$ is a convex set of dimension $\dim \sE(f'_v)$. Of course,
we may choose a bounded $U$. Then, 
picking $V$ as in Lemma~\ref{lemmafv},
we get $v+\epsilon U\subset V$,
for $\epsilon$ small enough, and applying~\eqref{neighbor},
we get 
\[
f(v+\epsilon u)=f(v)+f'_v(\epsilon u)=\lambda+v+\epsilon u
\qquad \forall u\in U \enspace,
\]
which shows that $v+\epsilon U\subset \sE(f)$. 
Therefore, $\dim \sE(f)\geq \dim \sE(f'_v)$.
To show the other inequality, we take
a set $U'\subset \sE(f)$ such that $\sfr_C(U')$
is a convex set of dimension $\dim \sE(f)$.
Then, by convexity of $\ec(f)$, for
all $0<\epsilon<1$, 
$\sfr_C((1-\epsilon) v+\epsilon U')\subset \ec(f)$.
We may choose a bounded $U'$. Then,
taking $V_\epsilon=\sfr_C^{-1}(\sfr_C((1-\epsilon)v+\epsilon U'))\subset \sE(f)$
and using the continuity of $\sfr_C^{-1}$ from $\ec(f)$ to $\sE(f)$
(which follows from the last assertion of Theorem~\ref{thmain}),
we get $V_\epsilon \subset V$ for $\epsilon$ small
enough, and using~\eqref{neighbor} again,
we get, for all $w\in V_\epsilon$,
\[
w-v= f( w)-\lambda-v= f(v)
+f'_v(w-v)-\lambda -v=f'_v(w-v) \enspace, 
\]
which shows that $V_\epsilon -v \subset \sE(f'_v)$.
Therefore, $\dim \sE(f'_v)\geq \dim \sE(f)$.
\end{proof}
We next give a new proof of the following result,
equivalent forms of which were proved by
Romanovsky, and Schweitzer and Federgruen.
\begin{theorem}[{Cf.~\cite[Th.~3]{roma73}, \cite[Th.~5.1 and Th.~5.3]{sch78}}]\label{th2}
The dimension of the eigenspace of a piecewise
affine convex monotone homogeneous map $f$
is equal to the number of critical classes,
$\sfm(f)$.
\end{theorem}
We shall reduce to the
special case of maps $g\geq I$ with eigenvalue $0$
(where $I$ denotes
the identity map). This case is simpler
because, when $g\geq I$,
\begin{equation}
\sE(g)=\set{x\in \R^n}{g(x)=x}=\set{x\in \R^n}{g(x)\leq x}=\sE^+(g)\enspace.
\label{eqgeqi}
\end{equation}
\begin{lemma}\label{lemmai}
If $g:\R^n\to\R^n$ is a piecewise affine convex monotone bihomogeneous
map such that $g\geq I$, then, $\sE(g)$
is a convex set of dimension $\sfm(g)$. 
\end{lemma}
Before proving Lemma~\ref{lemmai},
let us show that it implies Theorem~\ref{th2}.
Let $v$ denote any eigenvector of $f$,
and consider the maps $f'_v$ and $g=f'_v\vee I$, which both
are piecewise affine, convex, monotone and bihomogeneous.
Let $s=\sfm(f'_v)$, let $C_1,\ldots, C_s$ denote
the  critical classes of $f'_v$ (which are also
those of $f$), and $C=\nc(f'_v)=C_1\cup\ldots \cup C_s$.
It is immediate to check that $g$ also has $C_1,\ldots, C_s$
as  critical classes, and that is has
additional  critical classes 
of the form $\{i\}$, where $i$ is any element
of the complement $N=\{1,\ldots,n\}\setminus
C$. Thus, $\sfm(g)=s+|N|$. Applying
Lemma~\ref{lemmai} to $g$, we get that $\dim \sE(g)=s+|N|$.
Let $S$ denote any section of $f'_v$ (recall that sections were defined before
Corollary~\ref{cor-summarize}). Then, $S'=S\cup N$
is a section of $g$, and, by Corollary~\ref{cor-summarize},
$\sfr_{S'}(\sE(g))$ is a convex set of dimension
$s+|N|=|S'|$, which implies that $\sfr_S(\sE(g))$ is a convex
set of dimension $s$. 

Using 
Lemma~\ref{lem-projection}, we see that $(f'_v)^\omega$
sends $\sE(g)=\sE^+(f'_v)$ to $\sE(f'_v)$, and that 
$\sfr_S(\sE(g))$ coincides with  $\sfr_S(\sE(f'_v))$ ($S\subset C$).
 Therefore, $\sfr_S(\sE(f'_v))$ has dimension
$s$, which implies, by Corollary~\ref{cor-summarize},
that $\dim \sE(f'_v)= s=\sfm(f'_v)$.
By Lemma~\ref{reduc1}, we get 
$\dim \sE(f)=\dim \sE(f'_v)=\sfm(f'_v)=\sfm(f)$, which shows
Theorem~\ref{th2}. \qed

\begin{proof}[Proof of Lemma~\ref{lemmai}]
The convexity of $\sE(g)$ follows from~\eqref{eqgeqi}.
We shall prove 
that $\dim \sE(g)=\sfm(g)$
by induction on the dimension $n$. (The constructions
of the proof are illustrated in Example~\ref{ex0}
below).

If $n=1$, $g$ is the identity map, then,
$\dim \sE(g)=\sfm(g)=1$, and Lemma~\ref{lemmai}
holds trivially. Let us now assume that $n\geq 2$.
Denote by $F$ an invariant  
critical class of $g$, as defined in Lemma~\ref{univ}, 
let $N= \{1,\ldots,n\}\setminus F$, and
put $h=\widetilde{g_{NN}}$, where
$g_{NN}=\sfr_N \comp g\comp \sfi_N$,
with the notation of Section~\ref{secttang}.
Since $g\geq I$, we have $g_{NN}\geq I$, so $h\geq I$.

Corollary~\ref{cor-canon} allows us to write
for each $1\leq i\leq n$, $g_i(x)=\max_{p\in \sP_i}p\cdot x$,
where $\sP_i\subset \partial{g}_i(0)$ is a finite set of stochastic vectors.
Consider now, for all $i\in N$, the partition
$\sP_i=\sN_i\cup \sF_i $ where $\sN_i= \set{p\in \sP_i}{p_j=0,\forall j\in F}$,
and $\sF_i= \sP_i\setminus \sN_i$.
We have 
$h_i(y)=\max_{p\in \sN_i} p\cdot \sfi_N (y)$
for all $ i\in N, y\in \R^N $
(in particular, $\sN_i\neq \varnothing$,
for all $i\in N$, because $h\geq I$
takes finite values), and, for all $x\in \R^n$,
\begin{subequations}
\label{se-gsup}
\begin{eqnarray}
\label{se-gsup0}
g_i(x)&=&
h_i\comp\sfr_N(x) \vee \max_{p\in \sF_i} p\cdot x\qquad \forall i\in N,\\
g_i(x)&=& (g_{FF}\comp \sfr_F(x))_i \qquad 
\forall i\in F
\label{se-gsup1}
\enspace .
\end{eqnarray}
\end{subequations}
From~\eqref{e-disjointu}, we get
the disjoint union $\gc(g)=\gc(h)\cup \gc(g_{FF})$,
so $\sfm(g)=\sfm(h)+1$. By the induction assumption,
$\sE(h)$ is a convex set of dimension $\sfm(h)$, which implies that we can
find a bounded convex set $U\subset \sE(h)$ of dimension $\sfm(h)$.
We shall complete the elements of $U$ to get eigenvectors of $g$ as follows.
To each $u\in U$,
and $\lambda\in \R$, associate the vector $v(\lambda,u)\in \R^n$
such that $\sfr_N(v(\lambda,u))=u$, 
and $\sfr_F(v(\lambda,u))=\lambda$
(the constant vector).
We get from~\eqref{se-gsup1} and Lemma~\ref{lem-restriction},
\begin{equation}
g_i(v(u,\lambda)) = (g_{FF})_i(\lambda)=\lambda=
v(u,\lambda)_i
\qquad 
\forall i \in F 
\enspace ,
\label{eq-ok0}
\end{equation}
and we get from~\eqref{se-gsup0},
\begin{equation}
g_i(v(u,\lambda))= h_i(u ) \vee 
\max_{p\in \sF_i} p\cdot v(u,\lambda)=
u_i \vee 
\max_{p\in \sF_i}p\cdot v(u,\lambda)
\quad \forall i\in N 
 \enspace.
\label{eq-ok}
\end{equation}
Since for all $i\in N$ and $p\in \sF_i$, 
$p_j\neq 0$ for some $j\in F$,
we get $p\cdot  v(u,\lambda)= \sum_{j\in N} p_j u_j + \sum_{j\in F} p_j\lambda
\to -\infty$ when $\lambda\to -\infty$. Since $U$ is
bounded and the sets $\sF_i$ are finite, the limit is uniform
in $u\in U$ and $p\in \sF_i$. Therefore, we get from~\eqref{eq-ok}
that
\begin{equation}
g_i(v(u,\lambda))= u_i= (v(u,\lambda))_i 
\qquad 
\forall i\in N \enspace,
\label{eq-ok2}
\end{equation}
holds for all $u\in U$ and $\lambda\leq \lambda_0$, for some $\lambda_0\in \R$.
Combining~\eqref{eq-ok0} with~\eqref{eq-ok2}, we get
that $g(v(u,\lambda))=v(u,\lambda)$ for all $u\in U$ and $\lambda\leq \lambda_0$, which shows that $\sE(g)$ has dimension $\sfm(h)+1=
\sfm(g)$. 
\end{proof}
We obtain in passing some 
information on the restriction of $\sE(g)$ to
the set of invariant critical classes.
If $I\subset \{1,\ldots,n\}$, we denote by $1_I\in \R^n$
the vector such that $(1_I)_i=1$ if $i\in I$, and $(1_I)_i=0$
otherwise.
\begin{cor}
Let $g$ denote a convex piecewise affine 
monotone bihomogeneous map with invariant
critical classes $F_1,\ldots,F_r$.
Then, the restriction of $\sE(g)$
to $F_1\cup \cdots \cup F_r$
is equal to
$\R 1_{F_1} + \cdots + \R 1_{F_r}$.
\end{cor}
\begin{proof}
The reduction in the proof of Theorem~\ref{th2}
shows that we can assume that $g\geq I$,
so that we are in the situation of Lemma~\ref{lemmai}.
The corollary is obtained by a straightforward
variant of the induction argument of Lemma~\ref{lemmai},
in which $F$ is replaced by $F_1\cup \cdots \cup F_r$,
and $v(u,\lambda)$ is replaced by the vector $v(u,\lambda_1,\ldots,\lambda_r)$
defined by $\sfr_N(v(u,\lambda_1,\ldots,\lambda_r))=u$
and $\sfr_{F_i}(v(u,\lambda_1,\ldots,\lambda_r))=\lambda_i$ 
(the constant vector),
with $\lambda_1,\ldots,\lambda_r\in \R$, and $u\in U$.
The proof of Lemma~\ref{lemmai}
shows that for $\max_i \lambda_i$ close enough to $-\infty$,
$v(u,\lambda_1,\ldots,\lambda_r)$ is an eigenvector of $g$.
Then, for all $\lambda \in \R$, $v(\lambda+u, \lambda+\lambda_1,\ldots,
\lambda+\lambda_r)$ is an eigenvector of $g$, and
since $(\lambda+\lambda_1,\ldots,\lambda+\lambda_r)$ 
can take any value in $\R^r$,
this implies that the restriction of $\sE(g)$ to
$F_1\cup\cdots\cup F_r$ contains $\R 1_{F_1} + \cdots + \R 1_{F_r}$.
The other inclusion follows from the first assertion of Theorem~\ref{thmain}.
\end{proof}
\subsection{Examples}
\begin{example}\label{ex0}
Let $a_1,a_2,a_3>0$, and 
consider 
\[
f:\R^3\to\R^3, \qquad 
f(x)= \left(\begin{array}{c}
x_1 \vee (-a_1 + (x_2+x_3)/2)\\ 
x_2 \vee (-a_2 + (x_1+x_3)/2 )\\
x_3 \vee (-a_3 + (x_1+x_2)/2 )
\end{array}\right) \enspace .
\]
We have $f(0)=0$ and
$f'_0(x)= x$. Therefore, Theorem~\ref{th2} states
that $\sE(f)$ has dimension $3$.
Indeed, it is immediate to check that
$\sE(f)= \set{x}{x_i\geq -a_i + (x_k+x_j)/2, \forall 1\leq i\leq 3,
k\neq i, j\neq i}$, and this set has dimension $3$.
Since $\sE(f)$ is invariant by the translations
of vector $(\lambda,\lambda,\lambda)$, for all $\lambda\in \R$,
it is convenient to represent the projection 
of $\sE(f)$ on any plane orthogonal to the direction $(1,1,1)$,
which is as follows (we take $a_1=1/2, a_2=1, a_3=2$, the
point $(0,0,0)$ is represented by a bold point, the projection of
the eigenspace is the shaded region).
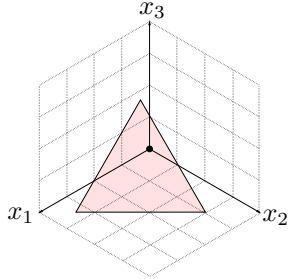
\begin{center}
\input fig6
\end{center}
By deforming this picture,
it should be obvious that in the limit case $a_1=a_2=0,a_3>0$,
the eigenspace $\sE(f)$, which looks as follows, 
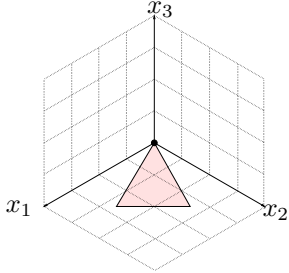
\begin{center}
\input fig7
\end{center}
has still dimension $3$. Let us
check this. We have in this case
\begin{equation}
f'_0(x)= \left(\begin{array}{c}
x_1 \vee (x_2+x_3)/2\\ 
x_2 \vee (x_1+x_3)/2 \\
x_3 
\end{array}\right) \enspace .
\label{ex11}
\end{equation}
We see that $f'_0$ has the three  critical classes
$\{1\}, \{2\}, \{3\}$, and, by Theorem~\ref{th2},
$\sE(f)$ has dimension $3$. It is instructive to illustrate
the proof of Lemma~\ref{lemmai} by the example
of $g=f'_0$. There is only one invariant critical class,  $F=\{3\}$.
We have $N=\{1,2\}$, and $h(y)=y$, for all $y\in \R^2$.
The representation~\eqref{ex11}
corresponds to $\sN_1=\{(1,0,0)\},\sF_1=(0,1/2,1/2)$,
$\sN_2=\{(0,1,0)\}$,
$\sF_2=(1/2,0,1/2)$. 
Any $u=(u_1,u_2)\in \R^2$
is an eigenvector of the identity map $h$.
The last step of the proof of the lemma
shows that for $\lambda$ negative enough, $v(u,\lambda)=(u_1,u_2,\lambda)$ 
is an eigenvector of $f'_0$. More precisely,
it is not difficult to see that 
the eigenspace of $f'_0$ looks as follows
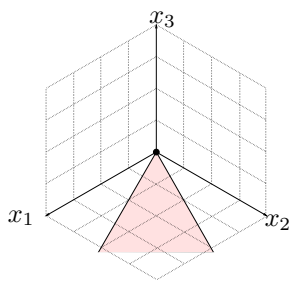
\begin{center}
\input fig8
\end{center}
where the shaded region now extends infinitely
in the direction $x_3\to -\infty$.
(In this example, $\sE(f'_0)$ is the {\em contingent
cone} to $\sE(f)$ at the point $0$.)
\end{example}
\begin{example}\label{excycli2}
We next present a variant of Example~\ref{ex0},
for which $f$ is of cyclicity $2$.
Let $a_1,a_2,a_3>0$, and 
\[
f:\R^3\to\R^3, \quad 
f(x)= \left(\begin{array}{c}
x_2 \vee (-a_1 + (x_1+x_3)/2)\\ 
x_1 \vee (-a_2 + (x_2+x_3)/2 )\\
x_3 \vee (-a_3 + (x_1+x_2)/2 )
\end{array}\right) \enspace .
\]
We have $f(0)=0$, and $\partial{f}(0)=\{f'(0)\}$ with
\[f'(0)=
\left(\begin{array}{ccc}
0 & 1& 0 \\
1 & 0& 0\\
0 & 0 & 1
\end{array}\right)\enspace,
\]
which gives the critical graph:
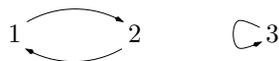
\begin{center}
\begin{tabular}[c]{c}
\input fig9
\end{tabular}
\end{center}
Therefore, $\gc(f)$ is the union
of two strongly connected components, $\sG_1$, and $\sG_2$,
with respective sets of nodes $\{1,2\}$
and $\{3\}$ and cyclicities $2$ and $1$,
which implies that $f$ has cyclicity $\sfc(f)=\lcm(2,1)=2$.
An immediate computation shows that 
\[
\sE(f)= \set{x\in \R^3}{ a_3 + x_3 \geq x_1=x_2\geq -2(a_1\wedge a_2)
+x_3} \enspace .
\]
This two dimensional convex
set is represented by the vertical bold segment on the following
figure:
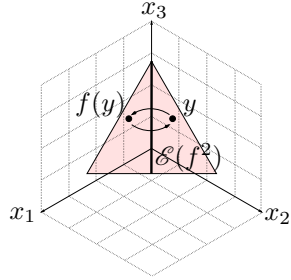
\begin{center}
\input fig10
\end{center}
Corollary~\ref{c-cycdim} says that the union of periodic
orbits of $f$ is equal to $\sE(f^2)$. A new
computation shows that $\sE(f^2)$ is equal to the triangular region
determined by the inequalities
\[
\begin{array}{lcl}
x_1& \geq&  -(a_1 \wedge a_2)+ (x_2+x_3)/2\enspace,\\
x_2 &\geq& -(a_1 \wedge a_2)+ (x_1+x_3)/2 \enspace,\\
x_3 &\geq&  -a_3 + (x_1+x_2)/2 \enspace ,
\end{array}
\]
which is depicted in the figure.
On this example, Theorem~\ref{th2} predicts that $\sE(f)$ is of dimension $2$,
and Proposition~\ref{orbitclass}, combined
with Theorem~\ref{th2}, predicts 
that $\sE(f^2)$ is of dimension $\sfc(\sG_1)+\sfc(\sG_2)=3$.
It is easy to see that $f$ acts on $\sE(f^2)$ as the reflection
with respect to the plane $x_1=x_2$. As an illustration, 
a periodic orbit $\{y,f(y)\}$ is depicted on the figure.
\end{example}
\begin{example}
Consider any convex monotone homogeneous
map $f:\R^n\to\R^n$, whose eigenspace  has dimension $n$.
We show here how, in that case,
Theorem~\ref{th2} allows us to characterize the set  $\sE(f)$.
We assume without loss of generality, that the eigenvalue is $0$.
Since the dimension of $\sE(f)$ is $n$,
the number of critical classes must
be equal to the dimension, $n$, hence $f\geq I$,
or equivalently, $f=I\vee g$, for some convex monotone homogeneous
map $g$. Then,
\begin{equation}
\sE(f)=\set{x\in \R^n}{x\geq Px-g^*(p)\;\forall P\in \dom g^*}
\enspace.
\label{e-carac-e}
\end{equation}
Half-spaces of the form $x_i \geq P_{i}x-g_i^*(p)$ involved
in this definition are special. Indeed, the later inequality
can be rewritten as
\begin{equation}
x_i \geq \sum_{j\neq i} \alpha_j x_j + \gamma \enspace,
\label{eq-infstable}
\end{equation}
with $0\leq \alpha_j$, $\gamma\in \R$,
and $\sum_{j\neq i} \alpha_j=1$. 
Using~\eqref{e-carac-e} and~\eqref{eq-infstable},
we get that a subset of $\R^n$ of non-empty interior is the eigenspace
of a convex monotone homogeneous map $f:\R^n\to \R^n$
if, and only if, it is the intersection of
half-spaces of the form~\eqref{eq-infstable}.
This implies, for instance, that the region
at the left hand side of the following
figure is the eigenspace
of a convex monotone homogeneous map:
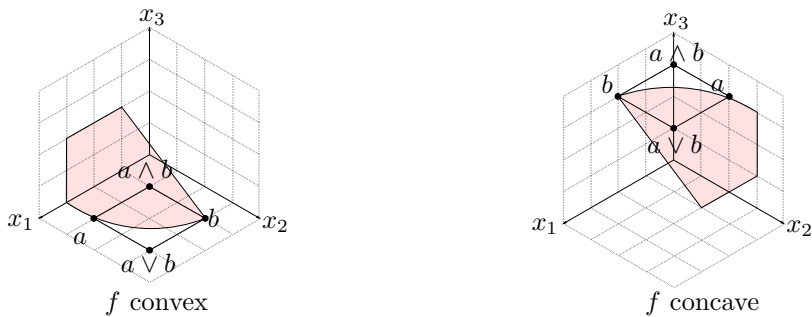
\begin{center}
\begin{tabular}{cc}
\input fig11&
\input fig12
\end{tabular}
\end{center}
The fact that $\sE(f)$ is an inf-subsemilattice
of $\R^n$ (but not a sup-subsemilattice)
is illustrated by the four points $a,b,a\vee b, a\wedge b$.
Eigenspaces of {\em concave} monotone homogeneous
maps have symmetric shapes, like the one at the right
hand side of the figure.
\end{example}
\begin{example}\label{ex-minmax}
The following example shows that for a non-convex
piecewise affine monotone homogeneous map $f$,
the dimension of the eigenspace is ill-defined.
Consider $f:\R^3\to \R^3$,
\begin{equation}
f(x)= \left(\begin{array}{c}
(x_1\wedge x_2 \wedge x_3) \vee (x_1\wedge (-1+x_2)\wedge
(1+x_3))\\
(2+x_1)\wedge x_2 \wedge (3+x_3)\\
(x_1\wedge x_2 \wedge x_3)
\end{array}\right) \enspace .
\label{ex111}
\end{equation}
Then, it is easy to check that $f=f^2$,
and that the eigenspace of $f$
is the following flag shaped set:
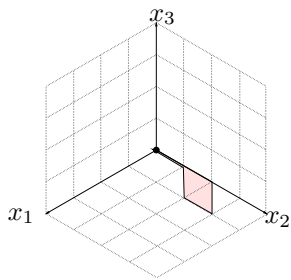
\begin{center}
\input fig13
\end{center}
Explicitly, $\sE(f)=K_1\cup K_2$,
where $K_1=\set{(\lambda,\lambda+t,\lambda)}{\lambda\in \R,0\leq t\leq 1}$,
and $K_2= \set{(\lambda,\lambda+t,\lambda+s)}{\lambda\in \R,1\leq t\leq 2,-1\leq s\leq 0}$,
which shows that the ``local dimension'' of the eigenspace near a point
$x\in \sE(f)$ is $2$ (if $x$ is in the relative
interior of $K_1$) or $3$ (if
$x$ is in the interior of $K_2$).
Although we do not need this here,
let us mention that there is a systematic technology to build
such examples, which originates from max-plus
algebra: $f$ is a projector on the
max-plus semimodule generated by the columns
of the matrix
\[
\left(\begin{array}{ccc}
0 & 0 & 0 \\
0 & 1 & 2\\
0 & -1 & -1
\end{array}\right) \enspace. 
\]
See~\cite{CGQ97a} for details.
\end{example}
\subsection{Computing $\gc(f)$}
\label{sect-algo}
To conclude this section, we give a polynomial
time algorithm to compute
\begin{equation}
\gf(\vex (\sQ_1 \times \cdots  \times \sQ_n))
\label{e-tcomp}
\end{equation}
given  finite sets of stochastic vectors  $\sQ_1,\ldots,\sQ_n$.
This algorithm allows us, in particular, to compute the critical graph
of a piecewise affine convex monotone homogeneous
map $f$, provided that an
eigenvector $u\in \R^n$ of $f$ is known.
Indeed, the coordinates of $f$ are of the form
\begin{equation}
f_i(x)= \max_{p\in \sP_i}(p\cdot x - f^*_i(p)) \enspace,
\label{sp3}
\end{equation}
where the $\sP_i$ are finite sets of stochastic vectors,
and, setting $\sQ_i= \set{p\in \sP_i}{p\cdot u- f^*_i(p)=f_i(u)}$,
it follows from \cite[Th.~10.31]{rocqw} that 
$\partial f_i(u)=
\vex \sQ_i $, hence $\partial f(u)=
\vex (\sQ_1\times \cdots \times \sQ_n)$.
Since $\gc(f)=\gf(\partial f(u))$, the problem
of computing $\gc(f)$ reduces to that of computing~\eqref{e-tcomp}.

To write the algorithm,
it will be convenient to consider more generally
a finite family $\{\sQ_i\}_{i\in I}$, where $\sQ_i\subset \R^I$
is a finite set of {\em substochastic} vectors.
For any  rectangular set $\sQ$ of {\em substochastic} matrices, 
we define $\gf(\sQ)$ as the union of the graphs of the matrices $P_{FF}$,
where $P\in \sQ$, $F$ is a final class of $P$, and $P_{FF}$ is a stochastic matrix
(this definition is consistent with
the one of $\gc(f)$ 
for a monotone subhomogeneous map $f$,
that we gave in~\S\ref{subsec-subhom}).

The algorithm can be specified as a recursive function,
with input $\{\sQ_i\}_{i\in I}$,
and output  $\gf(\vex \times_{i\in I} \sQ_i)$.
The function first builds, for all $i\in I$,
the subset $\sQ'_i\subset \sQ_i$
of row vectors with row sum $1$, together
with the graph $G=\sG(\times_{i\in I} \sQ'_i)=\cup_{P\in \times_{i\in I} \sQ'_i}
\sG(P)$.
If all the strongly connected components
of $G$ are trivial (we say that a strongly connected
component is \new{trivial} if it has only one node and no arcs),
the function returns the empty graph (with no nodes). Otherwise, we proceed
as follows. We denote by $F$ the union of final 
classes  of $G$, and put $N=I\setminus F$.
(The set $F$ is indeed the union of invariant final classes
of $f$ when $\vex \times_{i\in I} \sQ_i=\partial f(u)$, see Lemma~\ref{univ}.)
For $i\in N$, we define the sets $\sQ''_i\subset \R^N$ of row vectors
obtained by restricting to $N$ the vectors
$p\in \sQ'_i$ such that $p_j=0$, for all $j\in F$.
We denote by $G|_F$ the restriction of $G$ to $F$.  The identity
\begin{equation}\label{e-id}
\gf(\vex \times_{i\in I} \sQ_i) =
G|_F \cup  \gf(\vex \times_{i\in N} \sQ''_i)
\enspace,
\end{equation}
yields a recursive algorithm to compute
$\gf(\vex \times_{i\in I} \sQ_i)$. (The identity \eqref{e-id}
is similar to \eqref{e-disjointu}.)

Applying this algorithm
to the case of Example~\ref{exa1}, 
we get $\sQ_1=\{(1,0,0),(1/2,1/2,0),(1/2,0,1/2)\}$,
$\sQ_2=\{(0,1,0),(2/3,1/3,0)\}$,
$\sQ_3=\{(0,0,1)\}$,
$\sQ'_i=\sQ_i$,
$F=\{3\}$, $\sQ''_1=\{(1,0),(1/2,1/2)\}$,
$\sQ''_2=\{(0,1),(2/3,1/3)\}$, so that
\eqref{e-id} shows that $\gf(\vex (\sQ_1\times
\sQ_2\times \sQ_3))$ is the union 
of a loop at the node $3$, and of the complete
graph with nodes $1,2$.

\section{Stochastic Control Interpretation}
\label{sec-stoch}
In this section, we briefly explain how the above
results can be applied
to stochastic control. This application
also makes the results more intuitive.
See for instance ~\cite{whittle}, or~\cite{lasserre}
for more background on stochastic control.

A \new{Markov control model} with state space $\{1,\ldots,n\}$ 
is a $4$-uple $(A,\; \{A_i\}_{1\leq i \leq n},\;
\{r_i\}_{1\leq i \leq n},\; \{P_i\}_{1\leq i \leq n})$, where:
$A$ is a set, called \new{action space}; for 
each state $1\leq i\leq n$,  $A_i$ is a nonempty subset of $A$, whose elements
are interpreted as possible \new{actions};
$r_i$ is a map from $A_i$ to $\R$, the image $r_i^a$ of $a$ is
interpreted as an~\new{instantaneous reward}
received when action $a$ is performed in
state $i$; $P_i$ is a map from $A_i$ to the set  $\proba n$ of 
stochastic vectors, and the $j$-th entry, $P_{ij}^a$, 
of the image $P_i^a$ of $a$ is interpreted as the \new{transition probability}
from state $i$ to state $j$, when the action $a$ is performed. 
It will be convenient to assume that
$A$ is a topological space equipped with its Borel $\sigma$-algebra,
that the $A_i$ are Borel sets,
that the maps $r_i$ and $P_i$ are measurable, and that 
\begin{equation}
\mrm{for all $x\in \R^n$,}\quad 
 \sup_{a\in A_i} (r_i^a +P^a_i x)\;\;\;\mrm{is attained (and finite)}
\label{e-attained}
\end{equation}
(this is the case, in particular, if $A_i$ is compact, $r_i$ is
upper semi-continuous, and $P_i$ is continuous). 

The intuitive notion of {\em strategy}, i.e. of causal rule
telling which action to choose,
is captured by the following definitions.
An \new{history} is a sequence of states $(i_0,i_1,\ldots)$,
a \new{partial history} is a finite sequence $(i_0,\ldots,i_k)$.
A (randomized) \new{strategy} is a sequence
$\gamma=(\gamma_0,\gamma_1,\ldots)$ where $\gamma_k$ is
a map which to a partial history $(i_0,\ldots,i_k)$ associates
a probability measure $\gamma_k^{i_0,\ldots,i_k}$ on $A$ such that
$A_{i_k}$ has probability $1$, i.e. $\gamma_k^{i_0,\ldots,i_k}[A_{i_k}]=1$:
the action at time $k$, $a_k\in  A_{i_k}$,
will be chosen with probability $\gamma_k^{i_0,\ldots,i_k}$. We say
that $\gamma$ is \new{deterministic} if $\gamma_k^{i_0,\ldots,i_k}$ is 
a Dirac measure. We say that $\gamma$ is \new{Markovian} if
$\gamma_k^{i_0,\ldots,i_k}$ only depends of $k$ and $i_k$,
and that a Markovian $\gamma$ is \new{stationary}
if $\gamma_k$ is independent of $k$. Markovian
stationary policies are obtained by choosing,
for all $1\leq i\leq n$, a probability measure $\tau_i$ on $A$ such that $\tau_i[A_i]=1$, and taking $\gamma_k^{i_0,\ldots,i_{k}}= \tau_{i_k}$.
We denote by $\tau^\infty$ the (Markovian, stationary)
strategy $\gamma$ built in this way. 
For any strategy $\gamma$, and \new{initial state} $i_0$, we consider
a (state, action) stochastic process $(\xi_0,\alpha_0),(\xi_1,\alpha_1),\ldots$
with values in $\{1,\ldots,n\}\times A$,
such that $\xi_0=i_0$, the law of $\alpha_k$ knowing $\xi_0, \ldots, \xi_{k},
 \alpha_0\in A_{\xi_0}, \ldots, \alpha_{k-1}\in A_{\xi_{k-1}}$ is equal
to $\gamma_{k}^{\xi_0,\ldots,\xi_{k}}$, and the law of $\xi_k$ knowing
$\xi_0, \ldots, \xi_{k-1},
 \alpha_0\in A_{\xi_0}, \ldots, \alpha_{k-1}\in A_{\xi_{k-1}}$ 
is equal to $P^{\alpha_{k-1}}_{\xi_{k-1} \xi_k } $.
When $\gamma$ is Markovian and stationary,
$\xi_k$ simply becomes a time homogeneous Markov chain 
with initial state $i$ and transition matrix $P^{\tau}$,
\[
P^{\tau}_{ij}= \int_A P^a_{ij} d \tau_i(a) \enspace .
\]

The \new{ergodic control problem} consists in finding
a strategy $\gamma$ which maximizes, for all initial
states $1\leq i\leq n$, the \new{mean reward}
per time unit: 
\begin{subequations}
\begin{gather}
\mu_i^\gamma= \liminf_{N\to \infty} \frac{1}{N} \EE^{\gamma,i}
(r_{\xi_0}^{\alpha_0} + \cdots + r_{\xi_{N-1}}^{\alpha_{N-1}})
\enspace .
\end{gather}
We denote by 
\begin{equation}
\mu_i =\sup_{\gamma} \mu_i^\gamma 
\end{equation}\label{e-sub}
\end{subequations}
the optimum mean reward.

A closely related problem is to solve, for large $N$, the {\em horizon
$N$ problem}, which consists in maximizing
\begin{subequations}
\begin{gather}
v_i^\gamma(N)=  \EE^{\gamma,i}
(r_{\xi_0}^{\alpha_0} + \cdots + r_{\xi_{N-1}}^{\alpha_{N-1}} + \phi_{\xi_N})\enspace,
\end{gather}\label{e-sub2}
\end{subequations}
where the \new{final reward} $\phi$ is a map from $\{1,\ldots , n\}$ to $\R$,
or equivalently  an element of $\R^n$.
We set
\begin{gather} 
v_i(N) =\sup_{\gamma} v_i^\gamma(N) \enspace .
\end{gather}
The vector $v(N)\in \R^n$ is called the \new{value function}.
When $\gamma=\tau^\infty$, we shall simply write
$\mu^\tau$ and $v^\tau$, instead of 
$\mu^{\tau^\infty}$ and $v^{\tau^\infty}$,
respectively.

The study of both the ergodic control and finite horizon problems
relies on the \new{dynamic programming operator}, or {\em Hamiltonian},
$f:\R^n\to \R^n$,
\begin{equation}
f_i(x) = \sup_{a\in A_i} (r_i^a +P_i^a x) \enspace .
\label{e-mdp}
\end{equation}
It is obvious from~\eqref{e-mdp}
(and well known, see e.g.~\cite[Chap.~22, Th.~6.1]{whittle}), 
that the map $f$ is monotone, (additively) homogeneous,
and convex. Conversely, the Legendre-Fenchel duality
theorem shows that a convex monotone homogeneous
map $f$ can be written as
\begin{equation}
f_i(x)= \sup_{p\in \dom f_i^*} (p\cdot x - f_i^*(p))\enspace.
\label{e-lg}
\end{equation}
Since $\dom f_i^*$ is included in the set of stochastic
vectors (Proposition~\ref{prop-0}), and since $p\mapsto f_i^*(p)$
is lower semi-continuous, this is clearly of the form~\eqref{e-mdp}
(although $\dom f_i^*$ need not be compact, note that
property~\eqref{e-attained} is satisfied).
Moreover, when the action space $A$ is finite,
$\dom f_i^*=\vex A_i$ is a polyhedron,
and, by Proposition~\ref{prop-zero}, $f_i$ is
piecewise affine. 
Therefore:
\begin{prop}
The dynamic programming operators of Markov control models 
with state space $\{1,\ldots,n\}$
are exactly the convex, monotone homogeneous maps $\R^n\to \R^n$.
Moreover, Markov control models 
with finite action spaces correspond to piecewise
affine maps. 
\qed
\end{prop}
(The representation~\eqref{e-lg} provides
of canonical form for a Markov control models,
in which one can choose, when in state $i$, the 
transition probability from state $i$, $p\in \dom f_i^*$.)

The value function $v(N)$ can be computed
recursively via the dynamic programming equation
\begin{equation}
v(0)=\phi,\qquad v(N)=f(v(N-1))
\label{e-r}
\end{equation}
(see e.g.~\cite[Th.~3.2.1]{lasserre}).
Therefore, an eigenvector $u$, with associated
eigenvalue $\lambda$, yields
a \new{stationary solution} of the dynamic
programming equation~\eqref{e-r}, $v(N)= N\lambda + u$,
corresponding to the final reward $\phi=u$. Such 
stationary solutions are of economic interest.
Indeed, when $f(u)=\lambda+u$, 
we set 
\begin{equation}
\overline{A}_i=\set{a\in A_i}{\lambda + u_i = r_i^{a} + P^{a}_iu}
\label{e-ma}
\end{equation}
($\ov A_i\neq \emptyset$ thanks to~\eqref{e-attained})
and build  a Markovian stationary strategy $\tau^\infty$ 
by picking any probability measure $\tau$
such that $\tau[\overline{A}_i]=1$.
A standard result, that we shall not prove
(results of this kind can be found in ~\cite[Chap.~31]{whittle}
or~\cite[Th.~5.2.4]{lasserre}),
states that such a strategy is optimal both for the ergodic
control problem and for {\em all} the finite horizon
problems with final reward $\phi=u$, which means, loosely speaking,
that taking $\phi=u$ makes it possible
for the player to behave (optimally) in the
short term as he would in the long term 
(see~\cite{yakovenko} for more details on the
economic interpretation).

Another motivation
for describing the set of eigenvectors comes from~\eqref{e-r}:
Corollary~\ref{cor-asymp} shows that 
\[ v(N\sfc(f)) -N\sfc(f) \lambda \;\mrm{converges, when $N\to\infty$,}
\]
and the possible values of the limit are precisely the eigenvectors
of $f^{\sfc(f)}$.

Critical classes also have a stochastic control
interpretation. 
If $\gamma=\tau^\infty$ is a Markovian stationary strategy, 
we define the vector $r^\tau\in (\R\cup\{-\infty\})^n$ by 
\[
r^\tau_i = \int_A  r^a_i d\tau_i(a) \enspace .
\]
If $F$ is a final class of $P^\tau$, we denote by $m^{\tau}_F$ the
unique invariant measure of the $F\times F$ submatrix of $P^\tau$, and by 
$r^\tau_F\in \R^F$ the restriction of the vector $r^\tau$ to
$F$.
\begin{prop}\label{t-o}
Let us assume that the dynamic programming
operator $f$ has an eigenvector $u$, with associated eigenvalue $\lambda$.
Then, \begin{equation}
\lambda = \max_{\gamma\atop 1\leq i\leq n} \mu^\gamma_i 
=\max_{\tau^\infty\atop 1\leq i\leq n } \mu^\tau_i \quad\text{and}\quad
\max_{1\leq i\leq n } \mu^\tau_i
= \max_{F\;\mrm{final class of $P^\tau$}} m^{\tau}_F
r^{\tau}_F \enspace,\label{e-max}
\end{equation}
where the first two max are taken over
randomized strategies $\gamma$,
and randomized Markovian stationary strategies
$\tau^\infty$, respectively.
Moreover, if all the $A_i$ are compact,
and all the maps $r_i$ and $P_i$ are continuous, 
the elements of $\cc(f)$ are precisely the $F$ such that
$m^{\tau}_F r^{\tau}_F=\lambda$.
\end{prop}
Thus, the critical classes of $f$ are exactly
the maximal final classes
of randomized Markovian stationary strategies that are
optimal for the ergodic control problem.
(Considering {\em randomized} strategies is essential
for this equivalence to hold, even if the maximum
in~\eqref{e-max} is also attained by deterministic strategies.)
In particular, the \new{critical nodes}, i.e. the nodes of $\gc(f)$, 
are exactly the nodes which are visited infinitely often, almost surely,
by the trajectory of at least one optimal Markovian stationary policy.

\begin{proof}
The equalities in~\eqref{e-max} are standard results
of stochastic control (see for instance~\cite[Th.~5.2.4]{lasserre}
for the first two equalities; the last equality
follows from the ergodic theorem for reducible Markov chains, see
for instance~\cite[Ch.~31,\S6]{whittle}). 
We will only prove here the characterization of critical classes.

If $F\in \cc(f)$, $F$ is a final
class of a stochastic matrix $Q\in \partial{f}(u)$,
for some eigenvector $u$ of $f$. Subdifferentiating~\eqref{e-mdp}
at $x=u$, we get, thanks to the technical assumption on $A_i$, $r_i$ and $P_i$,
\[
\partial{f_i}(u) = \vex \set{P^a_i}{a \in \ov A_i}  \enspace,
\]
see~\cite[Th.~10.31]{rocqw}. Thus,
we can write $Q_i$ as a convex combination,
$Q_i=\sum_{a \in \ov A_i}  \alpha^a_i P^a_i$, where
the $\alpha^a_i$ are such that $\alpha^a_i\geq 0$,
$\sum_{a\in \ov A_i} \alpha^a_i=1$, and all but finitely many
$\alpha^a_i$ are zero.  Consider
now the measure $\tau_i=\sum_{a\in \ov A_i} \alpha^a_i\delta_{a}$
(the associated randomized Markov stationary strategy $\tau^\infty$
consists in playing the action $a$ with probability $\alpha^a_i$, when 
in state $i$, here $\delta_a$ denotes the Dirac probability
measure at $a$). Averaging the equalities $\lambda+ u_i 
=r^a_i + P_i^a u$ which hold for all $a\in \ov A_i$,
we get
\begin{equation}
\lambda + u = r^\tau + P^\tau u \enspace .
\label{e-kolm}
\end{equation}
Since $F$ is a final class of the matrix $P^\tau$,
restricting~\eqref{e-kolm} to $F$ and left
multiplying it by $m_F^\tau$, we get $\lambda +
m_F^\tau u_F= m_F^\tau r^\tau_F + m_F^\tau P_{FF}^\tau u_F$, and since
$m_F^\tau P^\tau_{FF}=m_F^\tau$, we get 
$\lambda =m_F^\tau r^\tau_F$.

Conversely, let us assume that $F$ is a final class of a matrix $P^\tau$
for some Markovian stationary strategy $\tau^\infty$,
and that $m^\tau_F r_F^\tau=\lambda$. For all $a\in A_i$, we set
\[
w^a_i=\lambda +u_i - r^a_i - P^a_i u  \geq 0  \enspace,
\]
and define the vector $w= (\int_A w^a_i d\tau_i(a))_{i\in F}$.
We have 
\begin{equation} \label{kolm2}
w=\lambda +u_F- r_F^\tau -P^\tau_{FF}u_F\enspace.
\end{equation}
Left multiplying~\eqref{kolm2} by $m_F^\tau$ and using again the fact that $m_F^\tau
 P^\tau_{FF}=m_F^\tau$,
together with  $m^\tau_F r_F^\tau=\lambda$,
we get $m_F^\tau w=0$. Since all the entries of $m_F^\tau$ 
are positive, and those of $w$ are nonnegative, 
we must have $w=0$. Hence
$ f_i(u) = r^\tau_i + P^\tau_i u$, for all $i\in F$, 
which implies that $P^\tau_i \in \partial{f_i}(u)$,
for all $i\in F$. Taking any element $Q\in\partial f(u)$ such that
$Q_i=P_i^\tau$ for $i\in F$, we obtain that $F$ is a final class of $Q$,
hence $F\in \cc(f)$.
\end{proof}

\end{document}